\documentclass[aap,MSNbibl,seceqn,citesort,dvips]{arximspdf}

\doi{10.1214/10-AAP753}
\volume{22}
\issue{1}
\pubyear{2012}
\firstpage{29}
\lastpage{69}

\makeatletter

\newtheorem{Lemma}{Lemma}[section]
\newtheorem{lemma}{Lemma}[section]
\newtheorem{Proposition}[Lemma]{Proposition}
\newtheorem{Theorem}[Lemma]{Theorem}
\newproclaim{remark}[Lemma]{Remark}
\newproclaim{remarks}{Remarks}
\newcommand{\A}{\mathcal{A}}
\newcommand{\tildmu}{\tilde{\mu}}
\newcommand{\ztss}{z_t^{(1,2)}}
\newcommand{\prob}{\mathbb{P}}
\newcommand{\PP}{\mathcal{P}}
\newcommand{\DA}{\mathcal{DA}}
\newcommand{\CC}{\mathcal{C}}
\newcommand{\tilBP}{\widetilde{\mathrm{BP}}}
\newcommand{\GG}{\mathcal{G}}
\newcommand{\NN}{\mathcal{N}}
\newcommand{\eps}{\varepsilon}
\newcommand{\Rbold}{\mathbb{R}}
\newcommand{\R}{\mathbb{R}}
\newcommand{\Zbold}{\mathbb{Z}}
\newcommand{\KK}{\mathcal{K}}
\newcommand{\expec}{\mathbb{E}}
\newcommand{\argmin}{\operatorname{arg\min}}

\newcommand{\e}{e}

\newcommand{\SWG}{\mathrm{SWG}}
\newcommand{\BP}{\mathrm{BP}}
\newcommand{\Op}{O_{\prob}}
\newcommand{\convp}{\stackrel{{\mathbb P}}{\longrightarrow}}
\newcommand{\convas}{\stackrel{\mathrm{a.s.}}{\longrightarrow}}
\newcommand{\convd}{\stackrel{d}{\longrightarrow}}

\newcommand{\vep}{\varepsilon}

\newcommand{\eqref}[1]{(\ref{#1})}

\makeatother

\begin{document}
\begin{frontmatter}

\title{Weak disorder asymptotics in the stochastic mean-field model of distance}
\runtitle{Weak disorder asymptotics}

\begin{aug}
\author[A]{\fnms{Shankar} \snm{Bhamidi}\corref{}\ead[label=e1]{bhamidi@email.unc.edu}\thanksref{titl1}}
and
\author[B]{\fnms{Remco} \snm{van der Hofstad}\ead[label=e2]{rhofstad@win.tue.nl}\thanksref{titl2}}
\runauthor{S. Bhamidi and R. van der Hofstad}
\thankstext{titl1}{Supported by NSF Grant DMS-07-04159, by PIMS and
NSERC Canada, a UNC Research Council and Junior faculty award.}
\thankstext{titl2}{Supported in part
by Netherlands
Organisation for Scientific Research (NWO).}
\affiliation{University of North Carolina and Eindhoven University of Technology}
\address[A]{Department of Statistics\\
\quad and Operations Research\\
University of North Carolina\\
304 Hanes Hall\\
Chapel Hill, North Carolina 27599\\
USA\\
\printead{e1}} %adresu isvedimo komanda gale!
\address[B]{Department of Mathematics\\\quad and
Computer Science\\
Eindhoven University of Technology\\
P.O.\ Box 513\\
5600 MB Eindhoven\\
The Netherlands\\
\printead{e2}}
\end{aug}

% HISTORY:
\received{\smonth{2} \syear{2010}}
\revised{\smonth{11} \syear{2011}}

% ABSTRACT
%
\begin{abstract}
In the recent past, there has been a concerted effort to develop
mathematical models for real-world networks and {to} analyze various
dynamics on these models. One particular problem of significant
importance is to understand the effect of random edge lengths or costs
on the geometry and flow transporting properties of the network. Two
different regimes are of great interest, the \textit{weak} disorder regime
where optimality of a path is determined by the sum of edge weights on
the path and the \textit{strong} disorder regime where optimality of a
path is determined by the maximal edge weight on the path. In the
context of the stochastic mean-field model of distance, we provide the
first mathematically tractable model of weak disorder and show that no
transition occurs at finite temperature. Indeed, we show that
for every finite temperature, the number of edges on the minimal
weight path (i.e., the hopcount) {is} $\Theta(\log{n})$ and
satisfies a central limit theorem with asymptotic means and variances
of order $\Theta(\log{n})$, with limiting constants expressible in
terms of the Malthusian rate of growth and the mean of the stable-age
distribution of {an} associated continuous-time branching process.
More precisely, we take independent and identically distributed edge
weights with distribution $E^s$ for some parameter $s>0$, where $E$ is
an exponential random variable with mean~1. Then the asymptotic mean
and variance of the central limit theorem for the hopcount are
$s\log{n}$ and $s^2 \log{n}$, respectively. We also find limiting
distributional asymptotics for the value of the minimal weight path in
terms of extreme value distributions and martingale limits of branching
processes.
\end{abstract}

% KEYWORDS
%

\end{frontmatter}

%s1 ###
\section{Introduction}
\label{sec-int}
The last few years have witnessed an explosion in empirical data
collected on various real-world networks, including transportation
networks like road and rail networks and data transmission networks
such as the Internet. This has stimulated an intense inter-disciplinary
effort to formulate mathematical network models to understand
their structure as well as the evolution of such real-world networks.
Rigorously analyzing properties of these models and deriving
asymptotics as the size of the network becomes large is currently an
active area of modern probability theory.

In many contexts, these models are used to describe
transportation networks, and understanding their flow carrying
properties is of paramount importance. Real-world networks are
described not only by their graph structure, which give us information
about the links between vertices in the network, but also by
their associated \emph{edge weights} that represent cost or time
required to traverse the edge. Similar questions form the core of one
of the fundamental problems in the modern theory of discrete
probability, namely \textit{first passage percolation}. In brief, one
starts with a finite network model $\KK_n$ (e.g., the
$[-n,n]^2$ box in the integer lattice $\Zbold^2$). Each edge $e$ is
given some random edge weight $l_e$, usually assumed to be nonnegative,
independent and identically distributed (i.i.d.) across edges. We shall
sometimes refer to $l_e$ as the \textit{length} or \textit{cost} of the
edge $e$. For any two vertices $u,v\in\KK_n$, and a path~$P$ between
the two vertices, the cost $f(P)$ of the path is some function of
the edge weights on the path (see the next section where we describe
two natural regimes). The optimal path $P_{\mathtt{opt}}(u,v)$ between the
two vertices is the path that minimizes this cost function amongst all
possible such paths. Now fix two vertices in $\KK_n$, for example, in
the case of the two-dimensional integer lattice, the origin and the
point $(n,0)$. One is then interested in deriving properties of the
optimal path between these two vertices, at least as the size of the
network tends to infinity.

In the modern applied context, two particular statistics of this
optimal path are of importance:
\begin{longlist}[(a)]
\item[(a)] $f(P_{\mathtt{opt}}(u,v))$: the cost of the optimal path. In many
situations, this gives the cost of transporting a unit of flow between
the two vertices.
\item[(b)] $H(P_{\mathtt{opt}}(u,v))$: the number of edges in the optimal
path. This represents the amount of time that a message takes in
getting between the two vertices. The mental picture one should have is
that the network is transporting flow between various vertices via
optimal paths, and \emph{delay}, that is, the amount of time that
a message takes in getting between vertices is the number of edges or
\textit{hops} on the optimal path. Thus, this quantity is often referred
to as the \textit{hopcount}.
\end{longlist}

%s1.1 ###
\subsection{Weak and strong disorder}
When modeling random disordered systems, two cost regimes for the cost
$f(P)$ of a path $P$ are of interest, the \textit{strong} disorder
and\vadjust{\goodbreak}
\textit{weak} disorder regime. Throughout the discussion below, we start
with a connected network $\KK_n$ on $n$ vertices, with each edge
assigned edge weight $l_e$. Fix two vertices denoted by $1$ and $2$
(say chosen uniformly at random amongst all vertices). We are
interested in properties of the optimal path between these two
vertices. Let $\PP_{12}$ denote the set of all paths between vertices
1 and 2.

\textit{Weak disorder regime}:
This is the conventional setup where, for any
path $\PP\in\mathcal{P}_{12}$, the cost of the path is
%
%e1.1 ###
\begin{equation}
f_{\mathtt{wk\mbox{-}dis}}(P) = \sum_{e\in P} l_e.
\end{equation}
The optimal path, denoted by $P_{\mathtt{wk\mbox{-}dis}}$, is defined by
%
%e1.2 ###
\begin{equation}
P_{\mathtt{wk\mbox{-}dis}} = \argmin\limits_{P\in\PP_{12}} f_{\mathtt{wk\mbox{-}dis}}(P).
\end{equation}
In our setup, the optimal path will always be \emph{unique}. We are
then interested in the cost and hopcount of this optimal path.

\textit{Strong disorder regime}:
Here, for any path $P\in\PP_{12}$, the
cost of the path is given by
%
%e1.3 ###
\begin{equation}
f_{\mathtt{st\mbox{-}dis}}(P) = \max_{e\in P} l_e.
\end{equation}
As before, the optimal path, denoted by $P_{\mathtt{st\mbox{-}dis}}$, is defined by
%
%e1.4 ###
\begin{equation}
P_{\mathtt{st\mbox{-}dis}} = \argmin\limits_{P\in\PP_{12}} f_{\mathtt{st\mbox{-}dis}}(P).
\end{equation}

From a statistical physics viewpoint, one is interested in
parametrizing the above problem via a real-valued parameter say $\beta
$, often called the ``inverse temperature'' of the system, such
that as $\beta\to\infty$, we get the strong disorder regime, while
for finite values of $\beta$, we have the weak disorder regime. One
interesting way of parameterizing the above problem is to consider the
original graph $\KK_n$ with some edge random variables $w_e$ and
consider the model~$\GG_n(\beta)$ where each edge is given weight
$l_e(\beta) = \exp(\beta w_e)$. The $\beta\to\infty$ regime then
corresponds to the strong disorder regime with edge weights~$w_e$, the
$\beta= 0$ regime corresponds to the graph distance
regime (where each edge has fixed weight $1$), while finite positive
values of $\beta$ are supposed to model the weak disorder regime and
are meant to interpolate between the graph distance regime and the
strong disorder regime. What is of paramount interest is to understand
if and when a transition occurs, namely given some model~$\KK_n$ of
network on $n$ vertices and edge distribution $w_e\sim F$,
for example,\ the uniform or exponential distribution, is there some finite
value of $\beta$ for which a transition occurs from the weak disorder
regime to the strong disorder regime, where the graph begins to behave
as in the strong disorder regime? What are the properties of the
optimal paths in various regimes, and how does the hopcount scale as a
function of $\beta$, at least in the $n\to\infty$ large network
limit? Although a number of studies have been carried out at the\vadjust{\goodbreak}
simulation level (see, e.g.,\ \cite{weak-strong-diso} and the
references therein) to understand such models of disorder in the
context of various random graph models resulting in fascinating
conjectures, there has been no rigorous effort carried out to derive
results in this context.

Our goal is to formulate a solvable model in this context and to
exhibit how such questions have deep connections to the stable-age
distribution theory of continuous-time branching processes as
formulated by Jagers and Nerman; see, for example,\ \cite{Jage75}. Without
further ado, let us dive into the formulation of the model in our context.

%s1.2 ###
\subsection{Model formulation}
\label{sec:not}
Let $\KK_n$ be the complete graph with vertex set $[n]\equiv\{1,
\ldots, n\}$ and edge set $E_n=\{ij\dvtx i,j\in[n], i\neq j\}$. Each edge
$e$ is given weight $l_e=(E_e)^s$ for some fixed $s>0$, where
$(E_e)_{e\in E_n}$
are i.i.d.\ exponential random variables with mean 1. The optimal path
between two vertices is the path that minimizes the sum of weights on
that path, as in the weak disorder regime. In the context of the above
discussion of strong and weak disorder, $s=0$ corresponds to the graph
distance, while $s=\infty$ corresponds to the strong disorder regime
with edge weights $E_e$,
the parameter $\beta>0$ above is equal to $s$ and the random variable
$(w_e)_{e\in E_n}$ equals $w_e=\log{(E_e)}$, which has a Gumbel
distribution. The advantage of this formulation is that it gives a
model that can be rigorously analyzed. The $s=1$ regime is one of the
most well-studied models in probabilistic combinatorial
optimization (see, e.g.,\ \cite{aldous-assignment,aldous2003objective,boll-rior-gamar-sudak,frieze1985value,hofs-flood,janson123})
and often goes under the name of ``stochastic mean-field model of
distance.'' For a fixed $s\in\Rbold^+$, we are interested in various
statistics of the optimal path, in particular, in the
asymptotics for the weight and hopcount of the optimal path as $n\to
\infty$.

To state the results, we shall need some definition. Let $(Y_j)_{j\geq
1}$ be i.i.d.\ exponential random variables with mean $1$. Define the
random variables $L_i$ by the equation
%
%e1.5 ###
\begin{equation}\label{eqn:bp}
L_i = \Biggl(\sum_{j=1}^i Y_j \Biggr)^s.
\end{equation}
Let $\PP$ be the above point process, that is,
%
%e1.6 ###
\begin{equation} \PP= (L_1, L_2,\ldots). \label{eqn:pp-def}
\end{equation}
While the parameter $s$ plays an important role in our analysis, for
the sake of simplicity, we shall omit it from the notation. The reader
should keep in mind that all the important constructions that arise in
the analysis and in the description of our results, such as
the point process above, depend on this parameter.
Now consider the continuous-time branching process (CTBP) where
at time $t=0$ we start with one vertex (called the root or the original
ancestor), each vertex $v$ lives forever, and has an offspring
distribution \mbox{$\PP_v\sim\PP$} as in \eqref{eqn:pp-def} independently
of every other vertex. Let $(\BP_t)_{t\geq0}$ denote
the CTBP with the above offspring distribution. The general
theory of branching processes (see, e.g.,\ \cite{Jage75}) implies that
there exists a\vadjust{\goodbreak} constant $\lambda= \lambda(s)$, called the \emph
{Malthusian rate of growth}, that determines the rate of
explosive growth of this model. In particular, if $z_t = |\BP
_t|$ denotes the number of individuals born by time $t$, then there
exists a strictly positive random variable~$W$ such that
%
%e1.7 ###
\begin{equation}
\label{eqn:lrv}
\e^{-\lambda t}z_t\convas W,
\end{equation}
where $\convas$ denotes convergence almost surely. The constant
$\lambda$ satisfies the equation
%
%e1.8 ###
\begin{equation}
\label{eq:malthus-solve}
\sum_{i=1}^\infty\expec(\e^{-\lambda L_i} ) = 1.
\end{equation}
In this case, an explicit computation (see Lemma \ref{lemma:malthus} below)
implies that
%
%e1.9 ###
\begin{equation} \label{lambda-s-comp} \lambda=\lambda(s) = \Gamma
(1+1/s)^s.
\end{equation}

Now let $W^{(1)}, W^{(2)}$ be i.i.d.\ with distribution $W$
where $W$ is as defined above in \eqref{eqn:lrv}. Define the Cox
process $\PP_{\mathtt{cox}}$ which, given $W^{(1)}$ and $W^{
(2)}$, is a Poisson process on $\Rbold$ with rate function given by
%
%e1.10 ###
\begin{equation}\label{eqn:cox-pr}
\gamma(x) = \frac{2\lambda}{s} W^{(1)}W^{(2)} \e
^{2\lambda x},\qquad  x\in\Rbold.
\end{equation}
%
%where $\kappa=\kappa(s)$ is the constant
% \eqn{
% \lbeq{kappa-def}
% \kappa=\kappa(s)=\frac{1}{s}\int_0^\infty\int_0^\infty(x_1+
%x_2)^{1/s-1} \lambda\e^{-\lambda x_1} \lambda\e^{-\lambda x_2} dx_1
%dx_2.
% }
Let $\Xi^{(1)}$ denote the first point of the point process $\PP
_{\mathtt{cox}}$.

%{\bf Definition of mean constant $\beta(s)$:} Another important
%constant that appears in the formulation of the main result is defined
%as follows: Let $\mu(t)$ be the expected number of births by time $t$
%namely
% \eqn{
% \lbeq{mu-s}
% \mu_s(t) = \sum_{i=1}^\infty\prob(L_i < t)
% =\sum_{i=1}^\infty\int_0^{t^{1/s}} e^{-u} \frac{u^{i-1}}{(i-1)!} dt
% =\int_0^{t^{1/s}} 1 du=t^{1/s}.
% }
%Then define $\beta(s)$ by the equation
% \begin{equation}
% \beta(s) = \int_0^\infty\lambda(s) e^{-\lambda(s) t} t \cdot
% \end{equation}
%so that, by \refeq{mu-s},
% \eqn{
% \beta(s) = \int_0^\infty\lambda(s) e^{-\lambda(s) t} t^{1+1/s}dt
% =\lambda(s)^{-(1+1/s)}\Gamma(2+1/s).
% }
%For later purposes, we shall also compute $\beta(s)\lambda(s)$, which
%by \refeq{lambda-s} equals
% \eqn{
% \lbeq{beta-lambda}
% \beta(s)\lambda(s)=\lambda(s)^{-1/s}\Gamma(2+1/s)=\frac{
% =1+1/s=\frac{s+1}{s}.
% }

%s1.3 ###
\subsection{Results}
\label{sec-res}
We are now in a position to state our results. Recall that we started
with the complete graph where each edge has distribution $l_e=E_e^s$,
where $(E_e)_{e\in E_n}$ are i.i.d.\ exponential random variables
having mean one. The first result identifies the limiting distribution
of the weight of the minimal weight path while the second result\vspace*{-2pt}
below identifies the asymptotics for the number of edges on the minimal
weight path. In the statement below, $\convd$ denotes convergence in
distribution.

\begin{Theorem}[(The weight of the shortest-weight path)]
\label{theo:wt-shwp}
Let $\CC=\CC(s)$ denote the cost of the optimal path between vertices
$1$ and $2$. Then, as $n\rightarrow\infty$,
%
%e1.11 ###
\begin{equation} n^s\CC- \frac{1}{\lambda}\log{n} \convd2\Xi
^{(1)},
\end{equation}
and
%
%e1.12 ###
\begin{equation} \label{Xi-distr} 2\Xi^{(1)}\stackrel{d}{=} \frac
{1}{\lambda} \bigl(G - \log{W^{ (1)}}-\log{W^{(2)}} -\log{(1/s)} \bigr),
\end{equation}
%
%ERRG...}
where $G$ is a standard Gumbel random variable independent of $W^{
(1)}$ and~$W^{(2)}$,
and $W^{(1)}$ and $W^{(2)}$ are two independent copies of the
random variable~$W$ appearing in \eqref{eqn:lrv}.
\end{Theorem}

\begin{Theorem}[(CLT for the hopcount)]
\label{theo:hopcount}
Let $H_n=H_n(s)$ denote the hopcount, that is, the number of edges
on the optimal path between\vadjust{\goodbreak} vertices $1$ and $2$. Then, as $n\to
\infty$,
%
%e1.13 ###
\begin{equation}
\frac{H_n- s\log{n}}{\sqrt{s^2\log{n}}} \convd Z,
\end{equation}
where $Z$ has a standard normal distribution.
\end{Theorem}

\begin{remarks*}
(a) Our proof shows that the convergence in Theorems
\ref{theo:wt-shwp} and~\ref{theo:hopcount} in fact occurs \emph
{jointly}, namely
%
%e1.14 ###
\begin{equation} \biggl( n^s\CC- \frac{1}{\lambda}\log{n},\frac{H_n-
s\log{n}}{\sqrt{s^2\log{n}}} \biggr) \convd\bigl(2\Xi^{(1)},Z\bigr),
\end{equation}
where the limiting random variables $\Xi^{(1)},Z$
are \emph{independent.}

(b) Not much is known about the random variable $W$ in \eqref{eqn:lrv}.
Indeed, the branching property can be used in order to show that it satisfies
the distributional relation
%
%e1.15 ###
\begin{equation} \label{W-rec} W\stackrel{d}{=} \sum_{i=1}^{\infty}
\e^{-\lambda L_i} W_i,
\end{equation}
where $(W_i)_{i\geq1}$ is an i.i.d.\ sequence of random variables
with the same distribution as $W$ independent of
$(L_i)_{i\geq1}$, and where $L_i$ is defined in \eqref{eqn:bp}. Using~\eqref{W-rec} and properties of functionals of Poisson
processes, one can show that the function $\phi(u)= \expec(\e
^{-uW})$, defined for $u\in\Rbold^+$, is the unique function
satisfying the functional relationship
%
%e1.16 ###
\begin{equation} \phi(u) = \exp\biggl(\int_0^\infty{ [\phi(u\e
^{-\lambda x^s} ) -1]}\,dx \biggr),\qquad  \phi(0)=1.
\end{equation}
%
%side. Please check!}
When $s=1$, then one can see this way that $W$ is an exponential
random variable with rate 1, but for other values of $s$, we
have no explicit form of~$W$.

(c) The distributional equivalence given by \eqref{Xi-distr} is
proved in Lemma \ref{lemma:Xi-distr} below.
\end{remarks*}

%s1.4 ###
\subsection{Discussion}
\label{sec-dis}
In this section, we discuss the relevance of our results and how they
relate to existing literature as well as various conjectures from
statistical physics. The standing assumption in this discussion is that
optimal paths are uniquely defined.

\subsubsection*{First vs.\ second order results}
First order results\vspace*{-2pt} (in our context showing, for example, that $H_n/s\log
{n} \convp1$, where $\convp$ denotes convergence in probability) are
much easier to prove
than the detailed convergence in distribution proved in Theorems \ref
{theo:wt-shwp} and \ref{theo:hopcount}. One of the reasons for the
length of this paper is that proving second order distributional
convergence results in these sorts of problems proves to be much more
difficult. Further, while in previous studies (e.g.,\ \cite
{BhaHofHoo09b} for various random network models) the hopcount
satisfied a central limit theorem (CLT) with \textit{matching}
means and variances, Theorem \ref{theo:hopcount} is novel in the sense
that it says that, for\vadjust{\goodbreak} large $n$, the hopcount has an approximate
normal distribution with mean $s\log{n}$ and variance $s^2\log{n}$.
Theorems such as Theorem \ref{theo:wt-shwp} for the actual cost of the
minimal weight path have been proven in a number of contexts (see,
e.g.,\ \cite{BhaHofHoo09b,hofs-flood,janson123}), but often prove
quite tricky to handle due to the fact that we only recenter the random
variables and do not divide by a normalizing factor going to $\infty$.
Thus, one needs to be extremely careful in analyzing the contribution
of various factors as $n\to\infty$. See, for example,\ \cite{BhaHofHoo09b}
to see the various factors that could contribute to the limiting
distribution in the
context of exponential weights on a \emph{random graph}.\vspace*{-3pt}

\subsubsection*{Strong disorder regime and minimal spanning trees}
Under strong disorder, it is easy to check using any of the standard
greedy algorithms for constructing minimal spanning trees that the
number of edges in the optimal path between any two vertices in the
network has the same distribution as the number of edges between the
two vertices in the \emph{minimal spanning tree} (with edge weights
$l_e$). More precisely, the optimal path between two vertices in the
strong disorder regime is identical to the path between the two
vertices in the minimal spanning tree.

In the context of our model, under strong disorder (``the $s=\infty$
regime'') what is known is that for the complete graph, the hopcount of
the optimal path $H(P_{\mathtt{st\mbox{-}dis}}) \sim\Theta_{\prob}(n^{1/3})$.
Here, for two sequences of random variables $(X_n)_{n\geq1}$ and
$(Y_n)_{n\geq1}$, we write $X_n=\Theta_{\prob}(Y_n)$ if
$X_n/Y_n$ and $Y_n/X_n$ are tight sequences of random variables.
This was first conjectured in \cite{weak-strong-diso} and recently
proven in \cite{addarioberry2006dms}. The above result in particular
shows that no transition occurs for finite values of $s$. It might be
interesting to analyze the above model when $s=s_n$ is a function of
$n$ and see when the strong disorder regime emerges ($s_n\to\infty$
regime) or the graph distance type behavior is preserved ($s_n\to0$).
In our proofs, we have kept formulas as explicit as
possible in order to be able to use them later on to study the strong
disorder case or the
graph distance limit. Let us now heuristically discuss the strong
disorder regime.\vspace*{-3pt}

\subsubsection*{Heuristics for strong disorder}
We see that the hopcount obeys a CLT with asymptotic mean and variance
equal to
$s\log{n}$ and $s^2\log{n}$, respectively. It is reasonable to expect
that the
CLT with asymptotic mean and variance equal to $s_n\log{n}$ and
$s^2_n\log{n}$
remains valid when $s_n$ is not too large. However, when $s_n$ is quite large,
then we should be in a phase that is close to the minimal spanning tree,
for which the hopcount scales like $n^{1/3}$ and has variance of order $n^{2/3}$
(since it is not concentrated). It would be of great interest to see
until what value
of $s_n$ the CLT with parameters $s_n\log{n}$ and $s^2_n\log{n}$
remains valid.
By the above, we see that for this, $s_n$ cannot grow faster than $n^{1/3}$
for this to be true.
In analogy to the scaling for the diameter of the Erd\H{o}s--R\'{e}nyi
random graph with edge probability $p=(1+\vep_n)/n$, which has size
$\vep_n^{-1} \log{(\vep^3_nn)}$ as long as $\vep_n\gg n^{-1/3}$
\cite{RioWor08}, one may wonder whether the hopcount
scales in leading order as $s_n\log{(n/s_n^3)}$, as long as $s_n\ll n^{1/3}$,
and where $s_n$ plays a similar role as $1/\vep_n$.\vadjust{\goodbreak}

\subsubsection*{Our choice of edge weights} If we rescale our weights by
$n^s$, the (expected) number of link weights that are at most $x$, from
a given vertex is, in our model, equal to
%
%e1.17 ###
\begin{equation} \label{edge-weights} (n-1)[1-\e^{-x^{1/s}/n}]
\approx
x^{1/s}.
\end{equation}
Thus, our weights are
chosen such that the weights obey a power law close to 0. In Internet,
the link weights are prescribed by the Internet Service Providers
(ISPs). Around 2000, CISCO, one of the main manufacturers of Internet
routers, has recommended to use the link weights in OSPF, the
Internet's intradomain routing protocol, that are proportional to the
inverse of the capacity or bandwidth of the link. This recommendation
has been followed by a many ISPs in order to optimally provision and
manage their networks.

Assuming that the link weights equals the
inverse bandwidth or capacity, our scaling relation in (\ref
{edge-weights}) is equivalent to the statement that the (expected)
number of links from a given vertex with capacity at least $B$ is close
to $B^{-1/s}$ for $B$ large. Thus, there is a power-law relation for
the link capacities in our model, and $1/s$ is the power-law exponent
in this relationship. By varying $s$, we can obtain any power-law
exponent. Unfortunately, in Internet, measuring the link capacities is
a notoriously hard problem, and, as a result, precise measurements of
their empirical properties are not available. Thus, while our model may
appear reasonable, we have no way of empirically verifying it.

\subsubsection*{Other edge weights} Note that in our context the
distribution of edge weights is $F(x)=1-\exp(-x^{1/s}) \sim x^{1/s}$
for $x$ close to zero. One would expect that the results in the paper
carry over rather easily to
edge weights with distribution function $F$ for which
$F(x)=x^{1/s}(1+o(1))$ when $x\downarrow0$. When~$F(x)$ has entirely
different behavior at $x=0$, other
properties might arise. Indeed, in our current setting, we see that
with high probability
the shortest-weight path traverses only through edges of weights of
order $n^{-s}$,
which is the size of the minimum of $n$ i.i.d.\ random variables with
distribution
$E^s$, where $E$ is exponential with mean 1. Thus, the benefit of using
edges of such
small weight vastly outweigh the fact that the path thus become longer
[i.e., has
$\Theta_{\prob}(\log{n})$ edges]. Now, when $F(x)=\e
^{-x^{-a}}$ for some $a>0$, then
the minimum of $n$ such random variables is $(\log
{n})^{-1/a}(1+o_{\prob}(1))$, so that
the minimal weight edge in the complete graph equals
$2^{-1/a}(\log{n})^{-1/a}(1+o_{\prob}(1))$. Here, we write that
$o_{\prob}(b_n)$ to denote a random variable $X_n$ which
satisfies that $X_n/b_n\convp0$. In particular, when $a>1$, we
cannot expect the optimal path
to have length $\Theta_{\prob}(\log{n})$, as already the
immediate path between
vertices 1 and 2 has smaller weight than any path of length $\log{n}$.

Moreover, it is not hard to see that the minimal \emph{two-step path}
between vertices 1 and 2 has weight $2^{1+1/a}(\log
{n})^{-1/a}(1+o_{\prob}(1))$,
so that the hopcount is with high probability\vadjust{\goodbreak} at most $2^{1+2/a}$. Thus,
this simple argument proves that the hopcount is
\emph{tight} for all $a>0$ (as is the case for the CM with infinite
mean degrees \cite{BhaHofHoo09a}). In \cite{BhaHofHoo11}, this
setting is investigated in more detail, and it is shown that, for most
value of $a$, the hopcount converges in probability to a constant.
Thus, it is clear that weights with distribution function
$F(x)=\e^{-x^{-a}}$
belong to a \emph{different universality class} as compared to edge
weights~$E^s$,
where~$E$ is an exponential random variable and $s>0$. This leads us to the
following general program:

\begin{quote}
\textit{Identify the universality classes for the weights in
first passage percolation on the complete graph}.
\end{quote}

\subsubsection*{Extensions of our results to random graphs}
A significant amount of work, both at the nonrigorous (\cite
{weak-strong-diso,havlin06,havlin05,lidia-opt} and the references
therein) as well as at the rigorous level \cite
{vcg-random-shanky,BhaHofHoo09a,BhaHofHoo09b,hofs2,hofs1}, has been
devoted to first passage percolation on random network models.
What is now generally expected is that in a wide variety of network
models and general edge costs, under weak disorder the hopcount scales
as $\Theta(\log{n})$ and satisfies central limit theorems as in
Theorem \ref{theo:hopcount}. We hope that the ideas in this paper
can also be applied to first passage percolation problems on various
random graphs, such as the \emph{configuration
model} (CM) with any given prescribed degree distribution $(p_k)_{k\geq
0}$. In \cite{BhaHofHoo09b}, first passage percolation
with exponential weights was studied
on the CM with finite mean degrees, and it is proved that similar
results as on the complete
graph hold in this case. Indeed, the hopcount satisfies a~CLT with
asymptotically \emph{equal} mean
and variance equal to $\lambda\log{n}$, where $\lambda$ is some
parameter expressible in terms of the degree
distribution. We expect that when putting exponential weights raised to
the power $s$ on the edges
changes this behavior, and the means and variances will become \emph
{different} constants times
$\log{n}$. While the behavior in \cite{BhaHofHoo09b} is remarkably
universal,\vadjust{\goodbreak} we expect that
for weights equal to powers of exponentials, when the variance of the
degrees is infinite,
the asymptotic ratio of mean and variance will be $s$ as on the
complete graph,
while for \emph{finite} variances degrees, the ratio may be different.

We see that the behavior of first passage percolation on the complete graph
with weights $E^s$ (as studied in this paper)
gives rise to CLTs for the hopcount with means and variances of order
$\log{n}$,
while weights with distribution function $F(x)=\e^{-x^{-a}}$ give rise to
bounded hopcounts, as is the case for the graph distance when all weights
are equal to 1. It would be of great interest to extend such results to
\emph{random graphs}. In particular, it would be of interest to
determine when the hopcount satisfies a CLT with asymptotic mean and
variance proportional to $\log{n}$, and when the hopcount behaves in a
similar way as the graph distance as studied for the CM in \cite
{hofs2,hofs3,hofs1}. This leads us to the following question:

\begin{quote}
\textit{How do the universality classes of first passage percolation
on the configuration model relate to those on the complete graph?}
\end{quote}

%s1.5 ###
\subsection{Proof idea and overview of the paper}
\label{sec-pf-idea}
For the sake of notational convenience, we shall rescale each edge
length by a factor $(n-1)^s$, so that each edge has distribution
$(Y_e)^s$, where $Y_e$ are exponential random variables with mean
$n-1$. This does not change the optimal path while the cost of this
path is scaled up by $(n-1)^s$. For the remainder of the paper,
we shall think of the edge weights as lengths which thus induce a \emph
{random metric} on the complete graph, and shall often refer to the
optimal path between two vertices as the \emph{shortest path} between
them. We are interested in the optimal path between vertices $1$ and
$2$. Consider water percolating through the network started
simultaneously from two sources, vertices $1$ and $2$, at rate one.
Then the first time of collision between the two flow processes, namely
the first time the flow percolating from vertex 1 sees a vertex already
visited by the the flow percolating from vertex 2 (or vice versa) gives
the shortest path between the two vertices. Let $z^{n,(1)}_t$ and
$z^{n,(2)}_t$ denote the number of vertices seen by the flow
cluster by time $t$ for the flow emanating from vertices $1$ and~$2$,
respectively. For large $n$, the flow clusters look like
independent versions of the CTBPs as formulated in Section \ref
{sec:not}, at least until they collide. A~coupling is rigorously
formulated in Sections \ref{sec-exp-single-vertex} and \ref
{sec:simul-flow}. Further, they collide only when both clusters reach
size $\Theta_{\prob}(\sqrt{n})$. At a heuristic level, at any
time $t$, the rate of collision $\gamma_n(t)$ in a small interval
$[t,t+dt)$ should be
%
%e1.18 ###
\begin{equation} \gamma_n(t) \propto\biggl( \frac{z^{n,(1)}_t
z^{n,(2)}_t}{n} \biggr)\,dt.
\end{equation}
%
%The above form emphasizes the contribution of the two cases that could
%arise when collision happens namely:

%(A) There is a birth from the flow cluster emanating from $1$ which
%happens with rate $N_n^{(1)}(t) dt$ and the new vertex born is
%chosen from the flow from $2$, this happens with probability $N_n^{
%(2)}(t)/n$;

%(B) Viceversa.

Now we use the fact that for large $t$, $z_t^{n,(i)} \sim W^{
(i)} \e^{\lambda t},$
where $W^{(i)}$ is the limiting random variable for the associated
CTBP defined in
\eqref{eqn:lrv}, to see that
%
%e1.19 ###
\begin{equation} \label{lim-res-CTBP} \gamma_n(x) \propto\frac
{W^{(1)}W^{(2)} \e^{2\lambda x}}{n}.
\end{equation}
Thus, collisions happen at time $(2\lambda)^{-1}\log{n}\pm\Op(1)$,
where $\Op(b_n)$ denotes a~sequence of random variables $(X_n)_{n\geq
1}$ for which $|X_n|/b_n$ is a tight sequence. If we let
$T_{12}$ denote the collision time, then the length of the optimal path
equals $W_n = 2T_{12}$. The above argument gives asymptotics for the
collision time and hence the length of the optimal path.

For the hopcount, we shall use general branching processes arguments to
show that at large time $t$, if one is interested in the distribution
of the generations (in our context this gives the number of individuals
at various graph distances away from the root, namely the originating
vertices $1$ and~$2$), the contribution to the population comes from
generations $t/\beta(s)$ and the deviations are normally distributed
around this value. Here the constant $\beta(s)>0$ denotes the mean of
the stable-age distribution of the associated branching process.
Intuitively, the optimal path between vertices $1$ and $2$ as
constructed via the above simultaneous flow picture looks like the
following: Suppose\vadjust{\goodbreak} the connecting edge between the two clusters $(v_1,
v_2)$ arises due to the birth of a child to vertex $v_1$ in the flow
cluster of vertex $1$ and this child, $v_2$ has already been visited by
the flow from $2$. This happens at around time ${(2\lambda)^{-1}}\log
{n}\pm\Op(1)$. The hopcount $H_n$ of the optimal path is given by the equation
%
%e1.20 ###
\begin{equation} H_n = G_1+ G_2+1,
\end{equation}
where $G_1$ and $G_2$ are the generations of vertex $v_1$ and $v_2$ in
flow cluster $1$ and~$2$, respectively. Thus, understanding the
distribution amongst generations in the coupled branching processes
paves the way to understanding the hopcount. The remainder of this
paper involves the conversion of the above heuristic into a \emph
{rigorous argument.} The organization of rest of the paper is as follows:
\begin{itemize}
\item In Section \ref{sec:coupling}, we shall couple the simultaneous
flows from two vertices on $\KK_n$ with
CTBPs and show that the difference is negligible;
\item Section \ref{sec:opt-weight} shows that the above coupling
incorporated with technical results from CTBP theory give us
asymptotics for the recentered length of the optimal path, namely
Theorem \ref{theo:wt-shwp}.
\item Section \ref{sec:hopcount-kn} shows how the distribution of
individuals among different generations in the associated branching
process proves Theorem \ref{theo:hopcount}.
\item Finally, Section \ref{sec:ctbp} proves all the CTBP results we
need to carry out our analysis. This section is the most technical part
of the paper and the point of organizing the paper in this fashion is
to motivate the various results that are proved in Section \ref{sec:ctbp}.
\end{itemize}

%s2 ###
\section{Proofs}
\label{sec-pfs}
In this section, we prove our main results. Proofs of the necessary
CTBP results are deferred to Section \ref{sec:ctbp}.

%s2.1 ###
\subsection{Dominating graph flow by continuous-time branching processes}
\label{sec:coupling}
In this section, we describe a coupling between the flows started from
vertices~1 and 2 and their
corresponding independent CTBPs with offspring distribution given by
the point process in \eqref{eqn:pp-def}. We shall first start with the
flow started
from one vertex and then extend this to the simultaneous flow
from two vertices.

%s2.1.1 ###
\subsubsection{Expansion of the flow from a single vertex}
\label{sec-exp-single-vertex}

We start with some notation. Recall that $\KK_n$ denoted the random
disordered media represented by the complete graph where each\vspace*{1pt}
undirected edge $(i,j)$ has edge length $E_{ij}^s$ where $E_{ij}$ are
i.i.d.\ exponentially distributed with mean $n-1$ [alternatively, with
rate $1/(n-1)$]. These edge lengths make $\KK_n$ a metric space (with
random geodesics). Let the index set of $\KK_n$ be $[n]:=\{1,2,\ldots
,n\}$ and fix vertex $1$. Think of this vertex as an originator of flow
of some fluid which percolates through the whole network via the
geodesics at rate $1$. Let $i_1=1, i_2, \ldots\in[n]$ be the
vertices
in sequential order seen by the flow. For $t\geq0$, let~$\SWG_t^{
(1)}$ be the\vadjust{\goodbreak} shortest-weight graph between vertex $1$ and all the
vertices that can be reached from $1$ by shortest-weight paths of
length at most~$t$. More precisely, $\SWG_t^{(1)}$ consists of these
shortest-weight paths and the weights of all of the edges used for
them. Let $(E_j^i)_{i\geq1, j\geq1}$ be a doubly infinite
array of mean $1$ exponential random variables. Then, by the properties
of the extremes of $n-1$
i.i.d. exponential random variables, each with mean $n-1$, it is
easy to see that the neighbors of $1$ have distances from $1$
distributed as
%
%e2.1 ###
\begin{equation} \label{weights-a}
\PP_{n,1} = ( E^{1}_1 )^s, \biggl( E^{1}_1+\frac
{n-1}{n-2}E^{1}_2 \biggr)^s, \ldots.
\end{equation}
Similarly, the distribution of distances from vertex $i_k$ (the $k$th vertex reached by the flow from $1$) to vertices other than those
already seen by the flow, is distributed as
%
%e2.2 ###
\begin{equation} \label{weights-b} \PP_{n,k} = \biggl( \frac{n-1}{n-k}E^k_1
\biggr)^s, \biggl( \frac{n-1}{n-k}E^k_1+ \frac{n-1}{n-k-1}E^k_2 \biggr)^s, \ldots.
\end{equation}
Call the above the \textit{immediate neighborhood process} of vertex $k$.
Note that for each $k$, by the memoryless property of the exponential
distribution,
the identity of the end point of each edge in the above point process is
uniformly distributed among all $[n]\setminus\{i_1, i_2, \ldots i_k\}$
vertices which have
not been seen at the time when the flow hits vertex $i_k$. Our aim is
to couple this process with a CTBP with offspring distribution
given by the point process $\PP$ defined by
%
%e2.3 ###
\begin{equation} \PP= \{(E_1)^s, (E_1+ E_2)^s, (E_1+ E_2+ E_3)^s,
\ldots\}, \label{eqn:pp-ctbp}
\end{equation}
where $(E_i)_{i\geq1}$ are i.i.d.\ exponential rate 1 random
variables. Comparing~\eqref{eqn:pp-ctbp} with \eqref{weights-a} and
\eqref{weights-b}, we see that, intuitively, the $\SWG_t^{(1)}$
should be stochastically smaller than the corresponding CTBP driven by
offspring distribution~$\PP$. The reason is that when the flow starts,
then the number of edges it has to explore from vertex $1$ is $n-1$,
but as the $\SWG_t^{(1)}$ increases with time, the number
of edges originating from each new vertex is strictly smaller than
$n-1$ due to vertices already explored by the flow. Thus, the points
are being \emph{depleted.} We shall show that asymptotically for large
$n$, the difference is negligible. To do so, as the flow explores~$\KK
_n$, we shall enlarge the graph~$\KK_n$ with new artificial vertices
to compensate for
the fact that~$\SWG_t^{(1)}$ uses up vertices in $\KK_n$
and effectively counteracting the
\emph{depletion of points} effect. For this, we shall need the
following randomization ingredients:

\begin{longlist}[(iii)]
\item[(i)] The complete graph $\KK_n$ with random edge weights;
\item[(ii)] An infinite array of i.i.d.\ exponential random variables
$(E_{i,j})_{i\in[n], j\geq n+1}$ each with mean $n-1$;
\item[(iii)] An infinite sequence of independent branching process
$(\tilBP_i(\cdot))_{i\geq n+1}$, each driven by the offspring
distribution in \eqref{eqn:pp-ctbp}.
\end{longlist}

Before diving into the construction, we shall need the following simple
lemma which follows directly from the memoryless property of the
exponential distribution.\vadjust{\goodbreak}

\begin{lemma}[(Powers of exponential distributions)]
\label{lem:expo-cond}
\textup{(a)} Consider the random variable $E^s$ where $E$ has an exponential
distribution with mean $n-1$. Then, for any fixed $r> 0$, the
conditional distribution of $E^s\mid E^s> r$ equals that of $(\tilde
{E} + r^{1/s})^s$, where $\tilde{E}$ is an independent
random variable with exponential distribution with mean $n-1$.

\textup{(b)} Consider the surplus random variable $(E^s- r)\mid E^s>r$. This
random variable has the same
distribution as the first point of a Poisson point process with rate
%
%e2.4 ###
\begin{equation} \label{rate-powers}
\Lambda_r(x) = \frac{1}{s (n-1)}
(r+x)^{1/s-1},\qquad  x\geq0.\vspace*{-2pt}
\end{equation}
\end{lemma}

We shall use part (a) of Lemma \ref{lem:expo-cond} in the construction
of the coupling while we shall use part (b) in the proof of the
distributional result for the optimal weight. We start by proving Lemma
\ref{lem:expo-cond}.

\begin{pf}
Part (a) is immediate from the memoryless property of the
exponential random variable.
For part (b), we note that
%
%e2.5 ###
\begin{eqnarray}
\prob(E^s-r\geq x\mid E^s> r) &=&\prob\bigl(E\geq
(x+r)^{1/s}\mid E> r^{1/s}\bigr)\nonumber\\[-9pt]\\[-9pt]
&=&\e^{-[(x+r)^{1/s}-r^{1/s}]/(n-1)},\nonumber\vspace*{-2pt}
\end{eqnarray}
while the probability that a Poisson point process with rate
(\ref{rate-powers}) has no points before $x$ equals
%
%e2.6 ###
\begin{equation} \quad\e^{-\int_0^x \Lambda_r(y)\,dy} =\e^{-\int_0^x
\frac{1}{s (n-1)} (r+y)^{1/s-1}dy}=\e
^{-[(x+r)^{1/s}-r^{1/s}]/(n-1)}.\vspace*{-2pt}
\end{equation}
Thus, the first point of this Poisson point process has the same
distribution as
the conditional law $E^s-r \mid E^s> r$.
\end{pf}

\subsubsection*{Construction of the coupling} This proceeds via the following
constructions:
\begin{longlist}[(a)]
\item[(a)] \textit{Artificial inactive vertices}: Consider the flow traveling at
rate one from vertex $1$ on $\KK_n$. Let $z_{t}^{n,(1)}$ denote
the number of vertices in $\SWG_t^{(1)}$. To evoke branching
process terminology, we shall often refer to this as the number of
vertices \textit{born} in the flow cluster of $1$ by time $t$. For $1\leq
k\leq n$, we define the stopping times
%
%e2.7 ###
%e2.17 ###
\begin{equation}
T_k^n = \inf\{t\dvtx z_t^{n, (1)} = k\},\vspace*{-2pt}
\end{equation}
so that $T_1^n=0$. Now consider the flow from vertex $1$. For $k\geq
2$, when the $k$th vertex $i_k$ is discovered by the flow at
time $T_k^n$, create a new artificial vertex labeled by $n+k-1$.
Let $a(i_k)$ denote the vertex in $\SWG_{T_k^n}^{(1)}$ to which
vertex~$i_k$ is attached.\vspace*{-2pt} Then note that for all $i_j\neq a(i_k)\in
\SWG_{T_k^n}^{(1)}$, by Lemma~\ref{lem:expo-cond}(a) and,\vspace*{-2pt}
conditionally on $\SWG_t^{(1)}$, the edge lengths of edge
$(i_j,i_k)$ have distribution $([T_k^n-T_j^n]^{1/s}+E)^s$ where $E$ has an
exponential distribution with mean $n-1$.\vadjust{\goodbreak}

For the new artificial vertex $n+k-1$, we attach edge lengths
from each vertex $i_j\in\SWG_{T_k^n}^{(1)}$ of length
$([T_k^n-T_j^n]^{1/s}+ E_{j,n+k-1})^s$ where the $E_{j,n+k-1}$
are exponential random variables as described in the randomization
needed for the coupling,
and where we recall that $T_j^n$ denotes the time of discovery of
vertex $i_j$.
We shall think of the flow having reached a distance $t-T_j^n$ on this edge.
At the time of creation, we shall think of these artificial vertices as
\textit{inactive} as the flow has not yet reached this vertex. Think of
these vertices as part of the network and the flow trying to get to
them as well. Note that eventually the flow will reach these inactive
vertices as well. Whenever the flow reaches an inactive artificial
vertex, we shall think of this vertex becoming active, that is, it is
\emph{activated}.
Let $\A_t$ denote the set of active artificial vertices. For $k\geq
1$, let
%
%e2.8 ###
\begin{equation} \label{hitting-time-k-act}
T_k^{n,*} := \inf\{ t\dvtx |\A_t|= k \}
\end{equation}
be the time of activation of the $k$th artificial vertex. Note
that in this construction,
edges exist only between vertices in $[n]$ and artificial vertices, no
edges exist between
artificial vertices.

\item[(b)] \textit{Activation of artificial vertices}: Note that activation of
inactive vertices happens at times $T_k^{n,*}$ via an edge from a
vertex in $\SWG_{T_k^{n,*}} \subseteq[n]$ to an inactive artificial
vertex $d_k\geq n+1$. Suppose at this time the set of artificial
vertices (active and inactive) is $\{n+1, n+2, \ldots,
n+j({T_k^{n,*}})\}$. When $d_k$ is activated, the following
constructions are performed:

\begin{longlist}[(3)]
\item[(1)] Remove all the edges from vertices in $[n]$ to $d_k$ (other
than the one that the flow used to get to it);
\item[(2)] Create a new inactive artificial vertex
$n+j({T_k^{n,*}})+1$. Just as before, create edges between each vertex
$i\in[n]$ and vertex $n+j({T_k^{n,*}})+1$ with edge lengths
distributed as $([t-T_k^{n,*}]^{1/s} + E_{i, n+j({T_k^{n,*}})+1})^s$
and think of the flow as having already traveled $t-T_k^{n,*}$ on it;
\item[(3)] At this time, start a CTBP $\tilBP_k(\cdot)$ with $d_k$
as the ancestor. The vertices born in this branching process have no
relation to the flow on $\KK_n$ and associated inactive vertices. For
time $t> T_{k}^{n,*}$, we shall call all the vertices in $\tilBP
_k(t)$, other than $d_k$, the \emph{descendants} of vertex $d_k$ at
time $t$.
\end{longlist}

Let $\DA_t$ denote the set of all descendants of the associated CTBPs
of active artificial vertices at time $t$ and let
%
%e2.9 ###
\begin{equation}\label{eqn:ff-couple}
\BP_t^{(1)} = \SWG_t^{(1)} \stackrel{\cdot}{\cup}
\A_t \stackrel{\cdot}{\cup} \DA_t.
\end{equation}
Let $z_t^{(1)} = |\BP_t^{(1)}|$ denote the number of
vertices reached at time $t$. The following proposition identifies
properties of the above construction that will be crucial in our
analysis. We shall prove this proposition in detail since later
we shall use an almost identical proposition in the context of flow
from two vertices which we shall state without proof in Section \ref
{sec:simul-flow} below.
\end{longlist}

\begin{Proposition}[(Properties of the coupling)]\label{prop:coupling}
In the above construction, the following holds:
\begin{longlist}[(a)]
\item[(a)] The process $(\BP_t^{(1)})_{t\geq0}$ is a CTBP driven
by the point process $\PP$ in~\eqref{eqn:pp-ctbp}. The process $(\SWG
_t^{(1)})_{t\geq0}$ is the shortest weight graph process
of the flow emanating from vertex $1$. As is obvious from \eqref
{eqn:ff-couple}, there is stochastic domination in the sense that, for
all times $t\geq0$, a.s.
%
%e2.10 ###
\begin{equation} \SWG_t^{(1)} \subseteq\BP_t^{(1)}.
\end{equation}
In particular, $z_{t}^{n, {(1)}}=|\SWG_t^{(1)}|\leq
z_t^{(1)}=|\BP_t^{(1)}|$ for all $t$.

\item[(b)] Let $\lambda= \lambda(s)$ be the Malthusian rate of growth of
$\BP_t^{(1)}$ as defined by \eqref{lambda-s-comp}.
Then, given any $\eps>0$, there exists $C_\eps>0$ such that for times
$t_n = (2\lambda)^{-1}\log{n}- C_\eps$
%
%e2.11 ###
\begin{equation} \liminf_{n\to\infty}\prob(\A_{t_n} = \varnothing
) \geq1-\eps.
\end{equation}

\item[(c)] For any fixed $B\in\Rbold$, letting $t_n = (2\lambda)^{-1}\log
{n} +B$, the sequence of random variable $|\A_{t_n}| + |\DA_{t_n}|$
is a sequence of tight random variables. Since the processes
$(|\A_{t}| + |\DA_{t}|)_{t\geq0}$ are monotonically increasing in
$t$, \eqref{eqn:ff-couple} implies that
$\sup_{t\leq t_n} (z_t^{(1)} - z_t^{n, {(1)}})$ is tight
and, in particular, as $n\to\infty$,
%
%e2.12 ###
\begin{equation} \sup_{t\leq t_n} \biggl|\frac{z_t^{n, {(1)}}}{z_t^{(1)}}
-1 \biggr| \convp0.
\end{equation}
\end{longlist}
\end{Proposition}

Note that if $|\A_{t_n}| = 0,$ then $\SWG_t^{(1)} = \BP_t^{(1)}$
for all $t\leq t_n$, so that part~(b) yields that
there is little difference between the SWG and the CTBP up to time
$(2\lambda)^{-1}\log{n}- C_\eps$.

\begin{pf}
Part (a) is obvious from construction. To prove part (b), note
that by construction, if $z_t^{n,(1)} = k$, then the chance that
the next vertex is an artificial inactive vertex is exactly $k/n$.
Thus, if $z_{t_n}^{n,(1)} = k_n $ then
%
%e2.13 ###
\begin{equation} |\A_{t_n}|\stackrel{d}{=} \sum_{j=1}^{k_n} I_j,
\label{eqn:bern}
\end{equation}
where $I_j$ are independent Bernoulli$(j/n)$ random variables, that is,
$\prob(I_j=1)=1-\prob(I_j=0)=j/n$. Now to choose $C_\eps$, first
choose $C_\eps^*>0$ so small that $\exp(-C_\eps^*/2) > 1-\eps/2.$
Since $z_{t_n}^{n,(1)} \leq z_{t_n}^{(1)}$ and for the
process $(z_t^{(1)})_{t\geq0}$ the asymptotics \eqref
{eqn:lrv} hold, we can choose $C_\eps^*$ such that
%
%e2.14 ###
\begin{equation}\label{eqn:zn-less-sqrn}
\prob\bigl(z_{t_n}^{n,(1)} > C_\eps^* \sqrt{n}\bigr)<\eps/2.
\end{equation}
Then
\begin{eqnarray*}
\prob(|\A_{t_n}| > 0) &\leq&\prob\bigl(|\A_{t_n}| > 0,
z_{t_n}^{n,(1)} < C_\eps^* \sqrt{n} \bigr)+ \prob\bigl(z_{t_n}^{n,
(1)} > C_\eps^*
\sqrt{n}\bigr) \\
&\leq&\bigl(1-\exp(-C_\eps^*/2)\bigr)+ \eps/2 < \eps,
\end{eqnarray*}
where the second inequality follows using a Poisson approximation in
\eqref{eqn:bern} and \eqref{eqn:zn-less-sqrn}. This proves part (b).

Finally to prove part (c), we note the following:
\begin{itemize}
\item Using part (b), we choose $C_\eps$ so that with high probability
no artificial vertices have been activated by time $(2\lambda
)^{-1}\log{n} - C_\eps$;
\item Using \eqref{eqn:bern} and ideas similar to the above argument,
one can show that the number of active artificial vertices by time $t_n
= (2\lambda)^{-1}\log{n} + B$ can be stochastically dominated with
high probability by a Poisson random variable $X_{ B}$ with mean
$C(B)$ for some function $B\mapsto C(B)$.
\end{itemize}
These two observations together imply that with high probability
%
%e2.15 ###
\begin{equation} |\A_{t_n}| + |\DA_{t_n}| \preceq_{\mathrm{st}} \sum
_{j=1}^{X_{ B}} |\BP_j(B-C_\eps)|,
\end{equation}
where $\BP_j(\cdot)$ are independent CTBPs driven by $\PP$,
independent of $X_{ B}$ which is Poisson with mean $C(B)$ and
$\preceq_{\mathrm{st}}$ denotes stochastic domination. Since the right-hand
side is bounded a.s., this proves part (c).
\end{pf}

%s2.1.2 ###
\subsubsection{Simultaneous expansion and coupling}
\label{sec:simul-flow}

Let us now show how the above coupling can be extended to flow
originating from two vertices $1,2$ simultaneously. We shall couple the
flow to two independent CTBPs $(\BP_t^{(i)})_{i=1,2}$. All
the ingredients of randomness shall be the same as in the previous
section, namely, (i) the complete graph $\KK_n$ with random edge
lengths; (ii) the infinite array of exponential random variables
$(E_{i,j})_{i\in[n], j\geq n+1}$; and
(iii) the infinite sequence of independent CTBPs $(\tilBP_i)_{i\geq
1}$ driven by $\PP$. Think of flow now emanating from the
two sources $1,2$ simultaneously at rate one exploring the shortest
weight structure about the two sources. We shall stop the flow when
there is a \emph{collision}, that is, the flow from one vertex sees a
vertex seen by the flow from the other vertex. As before, we let $\SWG
_t^{(i)}$ denote the shortest weight graphs up to time $t$
explored by the flow from each source $i=1,2$ and let
%
%e2.16 ###
\begin{equation}
\SWG_t = \SWG_t^{(1)} \cup\SWG_t^{(2)}.
\end{equation}
Let $z_t^{n,{(i)}} = |\SWG_t^{(i)}|$ and $z_t^{n} =
z_t^{n,{(1)}}+ z_t^{n,{(2)}}$.
Now let $T_{k}^n$ denote the stopping time
\begin{equation}
T_k^n = \inf\{t\dvtx z_t^n = k\},
\end{equation}
so that now $T_2^n = 0$. Let the vertex discovered at time $T_k^n$ and
attached to one of the two flow\vadjust{\goodbreak} clusters be $i_k\in[n]$. We shall call
this the \emph{time of birth of the vertex $i_k$}. Extra care is
needed as subtle issues of double counting of edges may arise.

The construction proceeds as before via two ingredients:

(a) \textit{Artificial inactive vertices}:
By convention, we shall think of the edge between $1$ and $2$ to belong
to the flow from vertex $1$, so that vertex 2 immediately is one
neighbor short. To compensate for this shortage, at time~$0$, we
shall add a new artificial inactive vertex labeled by $n+1$. Compared
to the other artificial vertices this shall be special in the sense
that vertex~$1$ will not have an edge to this vertex (or the artificial
vertices that replace this vertex when the flow reaches this vertex).
At time $0$, attach an edge $(2,n+1)$ of random length $E_{2,n+1}^s$.\vspace*{1pt}
Now start the flow from the two sources on the vertex set $[n] \cup\{
n+1\}$. The flow percolates from these two sources on the (expanded)
network discovering new vertices, both actual vertices in $[n]$ as well
as artificial vertices. Let $\SWG_t^*$ denote this flow process with
$z_t^{n,*} = |\SWG_t^*|$ and let
%
%e2.18 ###
\begin{equation}
\tilde{T}_k^n = \inf\{t\dvtx z_t^{n,*} = k \}.
\end{equation}
Let $i_k$ denote the vertex discovered by the flow at time $\tilde
{T}_k^n$ (this vertex could either be an actual vertex in $[n]$ or an
artificial inactive vertex). Create a new artificial vertex labeled by
$n+k$. Now if $i_k$ is in $\SWG_t^{(2)}$ then remove all the
edges between $i_k$ and all the vertices in $\SWG^{(2)}_{\tilde
{T}^n_k}$ (namely \textit{real} vertices in the actual\vspace*{-2pt} graph $[n]$ which
are part of $\SWG_t^{(2)}$ that have already been explored by the
flow from $2$). (Do the exact opposite if $i_k\in\SWG_t^{(1)}$.)
The edges $(v,i_k)$ for $v\in\SWG^{(1)}_{\tilde{T}^n_k}$ are
quite special (see the beginning of Section \ref{sec:opt-weight}).
Call these the \textit{potential connecting} edges as these are the edges
through which collisions of the two flow clusters may happen.
Also perform the following constructions:
\begin{itemize}
\item If $i_k\neq n+1$ or any of the \textit{replacements} of $n+1$ (this
term is defined below), then attach edges between the artificial vertex
$n+k$ and all $i_j\in\SWG_{\tilde{T}_k^n}$ with edge lengths
$([\tilde{T}_k^n-T^n_{i_j}]^{1/s}+E_{i_j,n+k})^s$. The flow would have already
reached up to distance $(\tilde{T}_k^n-T_{i_j}^n)$ on this edge to this new vertex.
\item If $i_k = n+1$, then replace this by a new vertex $n+k$. This
vertex will be called a \textit{replacement} of the special artificial
vertex $n+1$. Also replacements of such replacements shall be called
replacements. Remove all edges from $i_j\in\SWG_{\tilde {T}_k^n}$ to
$i_k$ and add back edges from these vertices\vspace*{1pt}
\textit{excluding} vertex~$1$ to vertex $n+1$ with edge lengths
$([\tilde{T}_k^n-T^n_{i_j}]^{1/s}+E_{i_j,n+k})^s$. This\vspace*{-1pt} can be
understood by noting that the flow would have already reached up to
distance $(\tilde{T}_k^n-T_{i_j}^{n})$ on this edge to this new vertex.
\end{itemize}

Every new artificial vertex when it is born is \textit{inactive}. Whenever
the flow reaches an inactive artificial\vadjust{\goodbreak} vertex we shall think of this
vertex becoming active and belonging to the flow cluster from which
this artificial vertex was reached. Let $\A_t^{(i)}$ denote the
set of active artificial vertices corresponding to flow cluster $i=1,2$
at time $t$ and let $\A_t = \A_t^{(1)} \cup\A_t^{(2)}$ be
the set of artificial vertices. For $k\geq1$, we let
%
%e2.19 ###
\begin{equation} T_{k}^{n,*} := \inf\{t\dvtx|\A_t|= k\}
\end{equation}
be the time of activation of the $k$th artificial vertex. Note
that, as before, edges exist only between vertices in $[n]$ and
artificial vertices
in this construction, no edges exist between artificial vertices.

(b) \textit{Activation of artificial vertices}: The flow will
eventually reach inactive artificial vertices. When this happens say
that \textit{activation} happens. This happens at times $T_k^{n,*}$ via an
edge from a vertex in $\SWG_{T_k^{n,*}} \subseteq[n]$ to an inactive
artificial vertex $d_k\geq n+1$ from one of the two flow clusters. When
an artificial vertex gets activated, it belongs to the flow cluster
that activates it and so do all its \textit{descendants} (the notion of a
descendant is defined below). Suppose that at this time, the set of
artificial vertices (active and inactive) is $\{n+1, n+2,\ldots,
n+j({T_k^{n,*}})\}$. As described above, this inactive artificial
vertex is replaced by a new inactive artificial vertex with appropriate
edges and edge lengths.

Further, at this time, start the CTBP $\tilBP_k(\cdot)$ with $d_k$ as
the ancestor. The vertices born in this branching process have no
relation to the flow on~$\KK_n$ and associated inactive vertices. At
time $t>T_{k}^{n,*},$ we shall call all the vertices in $\tilBP
_k$ other than $d_k$ the \emph{descendants of vertex $d_k$}.

Let $\DA_t^{(i)}$ denote the set of all descendants of the
associated CTBPs of active artificial vertices at time $t$ in flow
cluster $i=1,2$ and define the processes
%
%e2.20 ###
\begin{equation} \label{eqn:ff-couple-simul}
\BP^{(i)}_t = \SWG_t^{(i)} \cup\A_t^{(i)} \cup
\DA_t^{(i)},\qquad  i=1,2.
\end{equation}
Let $z_t^{(i)} = |\BP^{(i)}_t|$. Finally, let $\BP_t= \BP
_t^{(1)}\cup\BP_t^{(2)}$ denote the full flow process. This
completes the construction of the coupling.

The following proposition collects the properties of our construction
that we shall need. It is analogous to Proposition \ref{prop:coupling}
and we shall not give a~proof. Recall that $T_{12}$ denotes the
collision time of the two flow processes.

\begin{Proposition}[(Properties of the coupling)]\label{prop:coupling-simul}
In the above construction, the following holds:
\begin{longlist}[(a)]
\item[(a)] The processes $(\BP_t^{(i)})_{t\geq0}$ are independent
CTBPs driven by the point process $\PP$ in \eqref{eqn:pp-ctbp}. The
process $(\SWG_t^{(i)})_{0\leq t\leq T_{12}}$ is the
shortest weight graph process of the flow emanating from vertex $i$
till the collision time. As is obvious from \eqref
{eqn:ff-couple-simul},\vadjust{\goodbreak} there is stochastic domination in the sense that
for all times $t\geq0$,
%
%e2.21 ###
\begin{equation} \SWG_t^{(i)} \subseteq\BP_t^{{(i)}}.
\end{equation}
In particular $z_{t}^{n,{(i)}} \leq z_t^{(i)}$ for all $t\geq
0$.

\item[(b)] Let $\lambda= \lambda(s)$ be the Malthusian rate of growth of
$\BP^{(i)}_t$ as defined in \eqref{lambda-s-comp}. Then, given
any $\eps>0$, there exists $C_\eps$ such that for times $t_n =
(2\lambda)^{-1}\log{n} - C_\eps$,
%
%e2.22 ###
\begin{equation}
\liminf_{n\to\infty} \prob\bigl(T_{12}> t_n, \bigl|\A
_{t_n}^{(1)}\bigr| = 0, \bigl|\A_{t_n}^{(2)}\bigr| = 0 \bigr) \geq1-\eps.
\end{equation}
Note that if $|\A_{t_n}^{(i)}|=0$ then $\SWG_t^{(i)} = \BP
_t^{(i)}$ for all $t\leq t_n$.

\item[(c)] For any fixed $B\in\Rbold$, let $t_n^* = (2\lambda)^{-1}\log{n}
+B$ and let $t_n = T_{12}\wedge t_n^*$. Then the sequence of random
variables $|\A_{t_n}^{(i)}| + |\DA_{t_n}^{(i)}|$ is a tight
sequence of random variables. Since the processes $(|\A_{t}| + |\DA
_{t}|)_{0\leq t\leq T_{12}}$ are monotonically increasing in $t$,
\eqref{eqn:ff-couple-simul} implies that $\sup_{t\leq t_n} (z_t^{
(i)} - z_t^{n, {(i)}})$ is tight, and, as $n\to\infty$,
%
%e2.23 ###
\begin{equation}
\sup_{t\leq t_n} \biggl|\frac{z_t^{n, {(i)}}}{z_t^{(i)}}-1 \biggr|
\convp0.\vspace*{3pt}
\end{equation}
\end{longlist}
\end{Proposition}

%s2.2 ###
\subsection{Analysis of the weight of the optimal path}
\label{sec:opt-weight}
Before proceeding to the main proposition in this section, we shall
derive an important property of the above construction. When a vertex,
say $v\in[n]$, is born into one of the flow process (to fix ideas say into
the flow cluster of vertex $1$) at some time~$t$, then note that the
edges it has at this time are
\begin{itemize}
\item edges to inactive artificial vertices.
\item edges to all vertices in $[n]\setminus\SWG_t$.
\end{itemize}
For any vertex $v\in\SWG^{(1)}_t$ and, for any vertex $u\in[n]$
born into the flow cluster originating from vertex $2$ at some later
time $s>t$, we say that the edge connecting $v$ to $u$ is assigned to
vertex $v$ and {\bf not} to $u$. Similarly, if vertex $u$ is born into
the flow cluster starting from $2$ before vertex $v$ which is born into
flow cluster from vertex $1$, then say that the edge $(u,v)$ is
assigned to vertex $u$. Now, for any time $t$ and any vertex $v\in\SWG
^{(1)}_t\subseteq[n]$, let $N_t(v)$ denote the
number of edges with end points in $\SWG_t^{(2)}$ which are
assigned to it. Similarly, for a vertex $i\in\SWG_t^{(2)}$,
$N_t(v)$ is the number of vertices in $\SWG_t^{(1)}$
assigned to it. Recall that our aim in sending the flow simultaneously
is to analyze the collision time, namely, the first time when an edge,
which we shall refer to as the \textit{connecting edge}, forms between the
two flow clusters.\vadjust{\goodbreak} For any given time $t$ and $v\in\SWG_t^{(i)},
i=1,2,$ define the (\emph{random}) set
%
%e2.24 ###
\begin{eqnarray} \label{NNti-def}
\NN_t(v) &=& \bigl\{u\in\SWG_t^{(3-i)}\dvtx
\mbox{edge $(u,v)$ assigned to $v$}\bigr\}\nonumber\\[-8pt]\\[-8pt]
&=&\bigl\{u\in\SWG_t^{(3-i)}\dvtx T_u>T_v\bigr\},\nonumber
\end{eqnarray}
where, from now on, we shall use $T_v$ to denote the time of birth of
vertex~$v$ into the flow process $(\SWG_t)_{t\geq0}$ and we
recall that $\SWG_t=\SWG_t^{(1)}\cup\SWG_t^{(2)}$.

The importance of these connecting edges is as follows: Fix some time~$t$ and vertices $i\in\SWG_t^{(1)}$ and $j\in\SWG_t^{
(2)}$ with $T_j> T_i$ so that the edge between them is assigned to
vertex $i$. Note that up till time $T_j$, the flow was proceeding on
the edge between them at rate $1$ from vertex $j$. Now at time~$T_j$
the flow has reached the edge from the opposite side (i.e., from vertex
$j$) and is proceeding through the edge from \textit{both} end points.
Thus, while the flow through all other non-\textit{potential connecting}
edges proceeds at rate $1$, the flow through this edge proceeds at rate
$2$. For any time $x+T_j$, and using Lemma \ref{lem:expo-cond}(b) with
$r=T_j-T_i$ and the fact that the flow now proceeds at rate $2$ and not
$1$, the intensity function for the formation of this edge at this time is
%
%e2.25 ###
\begin{equation} \lambda_{(i,j)}(x+T_j)=\frac{2}{s(n-1)} \bigl((T_j-T_i) +
2x\bigr)^{1/s-1},\qquad  x\geq0. \label{eqn:rate-ij}
\end{equation}
In particular, for $t\geq T_j$,
%
%e2.26 ###
\begin{eqnarray}
\lambda_{(i,j)}(t) &=& \frac{2}{s(n-1)} \bigl((T_j-T_i) +
2(t-T_j)\bigr)^{1/s-1}\nonumber\\[-8pt]\\[-8pt]
&=&\frac{2}{s(n-1)} \bigl((t-T_i) + (t-T_j)\bigr)^{1/s-1}.\nonumber
\end{eqnarray}
This fact leads to the following proposition.

\begin{Proposition}[(Collision time distribution)]
\label{prop:rate-finite-n}
If $T_{12}$ denotes the collision time, then with respect to the
filtration generated by the flow process, $T_{12}$ has the same
distribution as the first point of a Poisson point process with rate
function given by
%
%e2.27 ###
\begin{equation} \label{rate-collision} \lambda_n(t) = \frac{2}{s
(n-1)} \sum_{i\in\SWG_t^{(1)}} \sum_{j\in\SWG_t^{(2)}} ([t-T_j]
+ [t-T_i] )^{1/s-1}.
\end{equation}
\end{Proposition}

\begin{remark}[(Extension to other graphs)]
Note that a similar formula as the above remains valid for any finite
graph with i.i.d.\ $E_e^s$ edge weights where $E_e$ are exponential
random variables, where the sum over $e=(i,j)$ is restricted to
$(i,j)\in E_n$, that is, the sum is only taken over the edges of the graph.
This can be used to analyze more general random graph models.
\end{remark}

\begin{pf*}{Proof of Proposition \ref{prop:rate-finite-n}}
Using \eqref{eqn:rate-ij}, Lemma \ref{lem:expo-cond} and the
fact that for a finite number of independent Poisson point processes,
the first point to occur in any of these processes has the same
distribution as the first point in Poisson point process with rate
given by the sum of rates of the corresponding point processes, we have that
%
%e2.28 ###
\begin{eqnarray}
\lambda_n(t) &=& 2 \sum_{i\in\SWG_t^{(1)}} \sum
_{j\in\NN_t(i)} \frac{([t-T_i] + [t-T_j])^{1/s -1}}{s (n-1)}\nonumber\\[-8pt]\\[-8pt]
&&{}+ 2
\sum_{i\in\SWG_t^{(2)}} \sum_{j\in\NN_t(i)} \frac{([t-T_i] +
[t-T_j])^{1/s -1}}{s (n-1)},\nonumber
\end{eqnarray}
where we recall that $\NN_t(i)$ denotes the set of vertices in the
other flow cluster assigned to $i$. Now note that for every pair of
vertices $(i,j), i\in\SWG_t^{(1)}, j\in\SWG_t^{(2)}$
either $i \in\NN_t(j)$ or vice versa and only one of these facts can
happen. Rearranging the above equation gives the result.
\end{pf*}

We call the sum appearing in (\ref{rate-collision}) a \emph
{two-vertex characteristic.}
In Section~\ref{sec-as-two-vertex} below,
we shall prove the following result concerning
the convergence of the two-vertex characteristic:

\begin{Theorem}[(Convergence of CTBP two-vertex characteristic)]\label{theo:denom-as-convg}
Consider two independent CTBPs $(\BP_t^{(i)})_{t\geq0},
i=1,2,$ as before.
Let $W^{(i)}, i=1,2$, be the almost sure limits of $\e^{-\lambda
t}z_t^{(i)}$. Then,
%
%e2.29 ###
\begin{equation} \e^{-2\lambda t}\sum_{i\in\BP_t^{(1)}} \sum_{j\in
\BP_t^{(2)}} ([t-T_j] + [t-T_i] )^{1/s-1}\convas\lambda W^{(1)}
W^{(2)},
\end{equation}
where $W^{(i)}$ are the a.s.\ limits of $\e^{-\lambda t}|\BP_t^{
(i)}|$ and are i.i.d.\ with the same distribution as $W$ in \eqref{eqn:lrv}.
%and $\kappa_1(s)$ is the constant
%e^{-\lambda(x_1+x_2)} dx_1 dx_2\]
\end{Theorem}

Now we are ready to prove Theorem \ref{theo:wt-shwp}.

\begin{pf*}{Completion of the proof of Theorem \ref{theo:wt-shwp}}
First, consider the rate function $\lambda_n(t)$ of the collision time given
in Proposition~\ref{prop:rate-finite-n}. By Proposition~\ref
{prop:coupling-simul}, in the summation arising in this rate function,
we can replace the terms $\SWG_t^{(i)}$ by $\BP_t^{(i)}$ as
the effect on the rate function is asymptotically negligible, where
$\BP_t^{(i)}$ are the independent CTBPs that have been coupled
with $\SWG_t^{(i)}$ to understand the optimal path on $\KK_n$.
Note that while the law of these CTBPs is independent of $n$, their
realizations intrinsically depend on~$n$,\vadjust{\goodbreak} since we have used the
randomization in $\KK_n$ to construct the CTBPs. We will indicate this
dependence by adding a subscript $n$. By \eqref{eqn:lrv},
%
%e2.30 ###
\begin{equation} \e^{-\lambda t}\bigl|\BP_t^{(i)}\bigr| \convas
W_n^{(i)},
\end{equation}
where $W_n^{(i)}$ are independent and identically distributed as
the limit variable in \eqref{eqn:lrv}.
% Recall that $z_t^{(i)} = |\BP_t^{(i)}|$. Now, the sum in the
%rate function above can be rearranged to give
% \eqn{
% \frac{1}{s (n-1)} \sum_{i\in\BP_t^{(1)}}
% \sum_{j\in\BP_t^{(2)}} ( [t-T_j] + [t-T_i] )^{1/s-1}= \frac{\e^{2
%(2)}\sum_{i\in\BP_t^{(1)}} \sum_{j\in\BP_t^{(2)}} \frac{
% }
%where $\lambda= \lambda(s)$ is the Malthusian rate of growth
%parameter and we define the \emph{two-vertex characteristic} $

%Theorem \ref{theo:denom-as-convg} shows that
% \eqn{
% \sum_{i\in\BP_t^{(1)}} \sum_{j\in\BP_t^{(2)}} \frac{
Now, Theorem \ref{theo:denom-as-convg} implies that for any fixed $B>
0$ %and any $x\in[-B,B]$, we have,
%
%e2.31 ###
\begin{equation} \sup_{x\in [-B,B]}\biggl|\lambda_n\bigl((2\lambda)^{-1}\log{n}
+x\bigr)- \frac{2\lambda}{s} W_n^{(1)} W_n^{(2)} \e^{2\lambda x} \biggr|\convp
0.
\end{equation}
Comparing the above with the definition of the Cox process in \eqref
{eqn:cox-pr} completes the proof subject to Theorem \ref
{theo:denom-as-convg}. Theorem \ref{theo:denom-as-convg} is proved in
Section~\ref{sec-as-two-vertex}.
\end{pf*}

%where $\kappa=\kappa(s)$ is the constant
% \eqn{\kappa(s)= \int_0^\infty\int_0^\infty
% (x_1+ x_2)^{1/s-1} \lambda\e^{-\lambda x_1} \lambda\e^{-\lambda
%x_2} dx_1 dx_2.
% }
%We note that we can identify $\kappa(s)$ as follows:
% \eqn{
% \lbeq{kappa-ident}
% \kappa(s)= \lambda^2 \int_0^\infty
% y^{1/s}\e^{-\lambda y}dy=\lambda^{-1/s+1} \Gamma(\frac{1}{s}+1)=
% }
%where we made use of \refeq{lambda-s-comp}.

For future reference, we define the \textit{two-vertex characteristics}
$\chi_{(i,j)}(t)$ by
%
%e2.32 ###
\begin{equation} \chi_{(i,j)}(t)= ([t-T_j] + [t-T_i] )^{1/s-1}. \label
{eqn:two-vert-charac}
\end{equation}

We shall now quickly prove the distributional equivalence \eqref{Xi-distr}.
\begin{Lemma}[(The limit of the shortest weight)]\label{lemma:Xi-distr}
The first point $\Xi^{(1)}$ of the Cox point process with rate
$\gamma(\cdot)$ as in \eqref{eqn:cox-pr} satisfies the
distributional equivalence in (\ref{Xi-distr}).
% \eqn{
% \Xi^{(1)} \stackrel{d}{=}\frac{1}{2\lambda} (G-\log{W^{(1)}}-
% where $G$ is a standard Gumbel random variable independent of $\{W^{
%(i)}:i=1,2\}$ which are independent and identically distributed as
%the limits of branching processes driven by the offspring distribution
%$\PP$ as in Equation \eqref{eqn:lrv}
\end{Lemma}

\begin{pf}
Since $\Xi^{(1)}$ is the first point of the Cox process
with rate function~$\gamma$ in \eqref{eqn:cox-pr},
we have for any fixed $y\in\Rbold$, conditional on $W^{(1)},
W^{(2)}$,
%
%e2.33 ###
\begin{eqnarray}
\prob\bigl(\Xi^{(1)}> y\mid W^{(1)}, W^{(2)}\bigr) &=& \exp\biggl(-\int
_{-\infty}^y \gamma(x) \,dx \biggr) \nonumber\\[-8pt]\\[-8pt]
&=&\exp\biggl(-\frac{1}{s}W^{(1)}W^{(2)} \e
^{2\lambda y} \biggr),\nonumber
\end{eqnarray}
so that
%
%e2.34 ###
\begin{eqnarray}
\prob\biggl(\Xi^{(1)}> x-\frac{1}{2\lambda}\log{\frac
{W^{(1)} W^{(2)}}{s}}\Bigm| W^{(1)} , W^{(2)} \biggr)&=&\exp(-\e^{2\lambda
x})\nonumber\\[-8pt]\\[-8pt]
&=&\prob\bigl(G/(2\lambda)> y\bigr),\nonumber
\end{eqnarray}
where $G$ has the standard Gumbel distribution. This proves the result.
\end{pf}

%s2.3 ###
\subsection{Hopcount analysis}
\label{sec:hopcount-kn}
As before, we let $T_{12}$ be the collision time between the two flow
clusters and suppose the collision happens via the formation of an edge
$(v_1, v_2)$ where $v_1 \in\SWG_{T_{12}}^{(1)}$ and $v_2 \in
\SWG_{T_{12}}^{(2)}$. For $i=1,2$, let $G_i$ denote the number of
edges on the path from vertex $i$ to $G_i$ so that the hopcount is
given by
%
%e2.35 ###
\begin{equation} H_n = G_1+ G_2+1.\vadjust{\goodbreak}
\end{equation}
To prove Theorem \ref{theo:hopcount} it suffices to show that, for
every fixed
$r,x, y \in\Rbold$ and writing $t_n=(2\lambda)^{-1} \log{n}$,
%
%e2.36 ###
\begin{eqnarray}  \label{aim-joint-conv}
&&\prob\bigl(T_{12}\leq t_n+ r, G_1
\leq\lambda s t_n + x s\sqrt{\lambda t_n}, G_2 \leq\lambda s t_n +
y s\sqrt{\lambda t_n} \bigr)\nonumber\\[-8pt]\\[-8pt]
&&\qquad  \to F_{12}(r)\Phi(x) \Phi(y),\nonumber
\end{eqnarray}
where $F_{12}(\cdot)$ is the distribution of the random variable
$\Xi^{(1)}$ appearing in Theorem~\ref{theo:wt-shwp} and $\Phi(\cdot)$ denotes
the standard normal distribution function.

For fixed time $t$ and $v\in\SWG_t^{(i)}, i=1,2$, let $G(v)$
denote the number of edges in the optimal path between $v$ and vertex
$i$ which started the flow. For any fixed $x\in\Rbold$, let
%
%e2.37 ###
\begin{equation} \SWG_t^{(i)}(x) = \bigl\{v\in\SWG_t^{(i)}\dvtx G(v)\leq
\lambda s t + xs\sqrt{\lambda t}\bigr\}.
\end{equation}
By Proposition \ref{prop:rate-finite-n} and properties of a finite
number of Poisson processes, we have, for any fixed $t$,
%
%e2.38 ###
\begin{eqnarray}\label{joint-conv}
&&\prob\bigl(T_{12} \in[t,t+dt), G_1 \leq\lambda st + x
s\sqrt{\lambda t}, G_2 \leq\lambda s t + y s\sqrt{\lambda t} \mid
\SWG_t\bigr)\nonumber\hspace*{-30pt}\\[-8pt]\\[-8pt]
&&\qquad  =\exp\biggl(-\int_0^t \lambda_n(w) \,dw \biggr) \lambda
_n(t) \frac{\sum_{i\in\SWG_t^{(1)}(x)} \sum_{j\in\SWG
_t^{(2)}(y)}\chi_{ij}(t)}{\sum_{i\in\SWG_t^{(1)}} \sum_{j\in\SWG
_t^{(2)}}\chi_{ij}(t)}\, dt,\nonumber\hspace*{-30pt}
\end{eqnarray}
where $\chi_{ij}(t)$ is the two-vertex characteristic defined in
\eqref{eqn:two-vert-charac}
and $\lambda_n(t)$ is the rate defined in (\ref{rate-collision}).
Thus, to complete the proof of (\ref{aim-joint-conv}),
it is enough to show the following theorem.

\begin{Theorem}[(CLT from two-vertex characteristic)]\label{theo:two-vertex-charac-ratio}
The two-vertex characteristic satisfies the asymptotics, for
$t\rightarrow\infty$,
%
%e2.39 ###
\begin{equation} \frac{\sum_{i\in\SWG_t^{(1)}(x)} \sum_{j\in\SWG
_t^{ (2)}(y)}\chi_{ij}(t)}{\sum_{i\in\SWG_t^{(1)}} \sum_{j\in
\SWG_t^{(2)}}\chi_{ij}(t)} \convp\Phi(x)\Phi(y).
\end{equation}
\end{Theorem}

Theorem \ref{theo:two-vertex-charac-ratio} is proved in Section \ref{sec:ctbp}
and completes the proof subject to Theorem~\ref{theo:two-vertex-charac-ratio}.
In fact, together with Theorem \ref{theo:denom-as-convg},
(\ref{joint-conv}) proves the joint convergence of the length of the
optimal path and the hopcount as remarked upon below Theorem~\ref
{theo:hopcount},
where the limits are \emph{independent}.

%s3 ###
\section{Continuous-time branching process theory}
\label{sec:ctbp}
In Sections~\ref{sec:opt-weight}--\ref{sec:hopcount-kn}, we have
reduced the proof of our main results to
the proof of Theorems~\ref{theo:denom-as-convg} and~\ref{theo:two-vertex-charac-ratio}.
In this section, we prove Theorems~\ref{theo:denom-as-convg} and~\ref{theo:two-vertex-charac-ratio}. This section is organized as
follows. In
Section~\ref{sec-CTBP-singleBP}, we investigate properties of our CTBP.
In Section~\ref{sec-one-vertex}, we investigate one-vertex
characteristics. In Section~\ref{sec-as-two-vertex}, we analyze
the two-vertex characteristic and prove Theorem~\ref{theo:denom-as-convg}.
In Section~\ref{sec-two-vertex}, we compute the mean and variance of
generation-weighted two-vertex characteristics, and
in Section~\ref{sec-pf-CLT}, we derive a CLT for the two-vertex characteristic
and complete our proof of Theorem~\ref{theo:two-vertex-charac-ratio}.

%s3.1 ###
\subsection{Intensities and limiting parameters for a single CTBP}
\label{sec-CTBP-singleBP}
We shall first state and prove various results that we shall require
regarding a \textit{single} branching process. Let $\BP$ be a
continuous-time branching process driven by the offspring point process
$\PP$ [i.e., the points given by $(L_1, L_2, \ldots)$ as in~\eqref
{eqn:bp}] and let $\mu$ denote the mean intensity measure of this
point process, that is,\looseness=-1
%
%e3.1 ###
\begin{equation} \mu[0,t] =\expec(\#\{i\dvtx L_i \leq t\}).
\end{equation}\looseness=0
Now,
%
%e3.2 ###
\begin{equation} \label{mu-s}
\qquad \mu[0,t] = \sum_{i=1}^\infty\prob(L_i
\leq t) =\sum_{i=1}^\infty\int_0^{t^{1/s}} e^{-u} \frac
{u^{i-1}}{(i-1)!} \,du =\int_0^{t^{1/s}} 1 \,du=t^{1/s}.
\end{equation}
Define the Malthusian rate of growth $\lambda= \lambda(s)$ as the
unique positive constant such that the measure
%
%e3.3 ###
\begin{equation} \nu(dt) = \e^{-\lambda t} \mu(dt) \label{eqn:nu}
\end{equation}
is a probability measure. A simple computation shows that this is
equivalent to \eqref{eq:malthus-solve}. The following lemma collects
some properties of this probability measure and the constant $\lambda$.

\begin{Lemma}[(Identification of limiting parameters CTBP)]
\label{lemma:malthus}
\begin{longlist}[(a)]
\item[(a)] The constant $\lambda=\lambda(s)$ is given by (\ref{lambda-s-comp}).

\item[(b)] The probability measure $\nu(dt)$ is a Gamma distribution with density
%
%e3.4 ###
\begin{equation} \label{f-form}
f(t) = \frac{\lambda^{1/s}}{\Gamma (1/s )} \e^{-\lambda t}t^{1/s-1}.
\end{equation}
%
%In particular, if $S_k=\sum_{i=1}^k X_i$ where $X_i$ are independent
%and identically distributed as $\nu$ then $S_k$ has density
% \eqn{
% \lbeq{fk-form}
% f_k(t) =\frac{\lambda^{k/s}}{\Gamma(\frac{k}{s} )} \e^{-\lambda t}t^{
% }
\item[(c)] Let $\beta_1$ and $\beta_2$ denote the mean and the standard
deviation of $\nu$. Then
%
%e3.5 ###
\begin{equation} \label{beta1-form}
\beta_1 = (s\lambda)^{-1}, \qquad \beta_2 = \bigl(\sqrt{s}\lambda\bigr)^{-1}.
\end{equation}
\item[(d)] Let $\mu^{*j}$ denote the $j$-fold convolution of the measure $\mu$. Then
%
%e3.6 ###
\begin{equation} \label{mu*j-form}
\mu^{*j}(du)=\frac{u^{j/s-1}\lambda^{j/s}\,du}{\Gamma(j/s)}.
\end{equation}
\end{longlist}
\end{Lemma}

\begin{pf}
To prove part (a), note that since the sum of $i$
independent exponential random variables follows the gamma distribution,
a simple computation gives that
\begin{eqnarray}\label{lapl-trans}
1&=&\sum_{i=1}^\infty\expec(\e^{-\lambda L_i} ) =\sum
_{i=1}^\infty\int_0^{\infty} \e^{-\lambda t^s} \e^{-t} \frac
{t^{i-1}}{(i-1)!} \,dt \nonumber\\
&=&\int_0^{\infty} \e^{-\lambda t^s} \e^{-t}
\sum_{i=1}^\infty\frac{t^{i-1}}{(i-1)!} \,dt
=\int
_0^{\infty} \e^{-\lambda t^s} \,dt\nonumber\\[-8pt]\\[-8pt]
&=&\lambda^{-1/s}\int_0^{\infty} \e
^{-t^s} \,dt
=\lambda^{-1/s} s^{-1}\int_0^{\infty} \e^{-v} v^{1/s-1}
\,dv\nonumber\\
&=&\lambda^{-1/s} \Gamma(1/s)/s=\lambda^{-1/s} \Gamma
(1+1/s) \nonumber
\end{eqnarray}
as required. Parts (b) and (c) are trivial. To prove part (d) note that,
by~(\ref{mu-s}) and \cite{GraRyz65}, equation 4.634, we have
%
%e3.7 ###
\begin{eqnarray}  \label{mu*j-comp} \mu^{*j}(du)&=&du\, s^{-j}\int
_{u_1+\cdots+u_j=u} u_1^{1/s-1}\cdots u_j^{1/s-1}\,du_1\cdots
du_j\nonumber\\
&=&\frac{u^{j/s-1} s^{-j}\Gamma(1/s)^jdu}{\Gamma(j/s)}=\frac
{u^{j/s-1} \Gamma(1+1/s)^jdu}{\Gamma(j/s)} \\
&=&\frac{u^{j/s-1} \lambda^{j/s}du}{\Gamma(j/s)}.\nonumber
\end{eqnarray}
\upqed
\end{pf}

%s3.2 ###
\subsection{Analysis of single-vertex characteristic}
\label{sec-one-vertex}
We first state a general theorem for \textit{single} vertex
characteristics of the CTBP. Consider a function $\chi\dvtx\Rbold
^+\to\Rbold^+$ which is continuous almost everywhere and
which (a) increases at most polynomially quickly at $\infty$;
and (b) is integrable with respect to the Lebesgue measure near zero.
Let us call such functions \emph{regular single-vertex
characteristics}. For the branching process $\BP_t$, call
%
%e3.8 ###
\begin{equation}\label{eqn:usual-charac}
z_t^\chi= \sum_{j\in\BP_t} \chi(t-T_j)
\end{equation}
the \emph{branching process counted according to characteristic} $\chi
$. Branching processes counted by characteristics are some of the
fundamental objects studied by Jagers and Nerman, see, for example, \cite
{JagNer84}. For example, taking $\chi(x) =1$, we obtain
$z_t^\chi= |\BP_t|$, the size of the branching process at time $t$.
In order to investigate the hopcount, we will need to analyze
not just branching processes counted according to characteristics as
above but also \emph{generation-weighted} characteristics. Given a
regular single vertex characteristic $\chi$ and any fixed $a\in\R$, define
%
%e3.9 ###
\begin{equation} z_t^\chi(a)=\sum_{j\in\BP_t} a^{G(j)}\chi(t-T_j),
\label{eqn:gen-wtd-charac}
\end{equation}
where, as before, $T_j$ denotes the time of birth of vertex $j$, while
$G(j)$ denotes the \emph{height} or \emph{generation}
of vertex $j$. Given any characteristic $\chi$,
when we write $z_t^\chi$ without the argument $a$, we imply the
branching process counted in the usual way as in
\eqref{eqn:usual-charac}, while when we have an argument $a$, namely~$z_t^\chi(a)$, we refer to the branching process counted by a
generation-weighted characteristic as in \eqref {eqn:gen-wtd-charac}.

The following proposition is adapted from the general theory of CTBPs,
see, for example,\ \cite{Jage75},\vadjust{\goodbreak} Theorem 5.2.2, for part (a)
(or see the nice treatment in \cite{JagNer84}, Theorem 3.4). We shall
give a complete proof since branching processes counted by
generation-weighted characteristics have not been previously analyzed.
These constructions shall be crucial for us in order to
prove the CLT for the hopcount.

\begin{Proposition}[(Mean and co-variances of one-vertex characteristics)]
\label{prop:form-1-charac-exp}
For regular deterministic single-vertex characteristics $\chi$,
\begin{longlist}[(a)]
\item[(a)] the expectation $m_t^\chi(a) = \expec[z_t^\chi(a)]$ satisfies
%
%e3.10 ###
\begin{equation} m_t^\chi(a)=\expec[z_t^\chi(a)] = \int_0^t \chi
(t-u)\sum_{j=0}^\infty a^j \mu^{*j}(du),
\end{equation}
%
% (b) the variance of $z_t^\chi(a)$ satisfies the equation
% \eqn{
% \operatorname{Var}(z_t^\chi(a)) = \int_0^t h_a(t-u) \tilde{\mu}_{a^2}(du),
% }
% where $t\mapsto h_a(t)$ is the function
% \eqn{
% h_a(t) = a^2\operatorname{Var} (\int_0^t m_{t-u}^\chi(a) \PP(du) ),
% }
% and we define
% \eqn{
% \lbeq{tilde-mu-a-def}
% \tilde{\mu}_{a}(du)=\sum_{j=0}^\infty a^{j}\mu^{*j}(du).
% }
\item[(b)] the covariances between $z_t^{\chi_1}(a_1)$ and $z_t^{\chi
_2}(a_2)$ satisfy
%
%e3.11 ###
\begin{equation} \label{Cov-form}
\operatorname{Cov}\bigl(z_t^{\chi_1}(a_1),
z_t^{\chi_2}(a_2)\bigr) = \int_0^t h_{a_1,a_2}(t-u)\tilde{\mu}_{a_1a_2}(du),
\end{equation}
where $v\mapsto h_{a_1,a_2}(v)$ is the function
%
%e3.12 ###
\begin{equation} \label{ha-cov-def} h_{a_1,a_2}(v) =\frac
{a_1a_2}{s}\int_0^{v} u^{1/s-1}m_{v-u}^{\chi_1}(a_1)m_{v-u}^{\chi
_2}(a_2)\,du,
\end{equation}
and we define the generation-weighted intensity measure $\tilde{\mu
}_a$ by
%
%e3.13 ###
\begin{equation} \label{tilde-mu-a-def} \tilde{\mu}_{a}(du)=\sum
_{j=0}^\infty a^{j}\mu^{*j}(du).
\end{equation}
\end{longlist}
\end{Proposition}

\begin{pf}
The proof of part (a) follows the same strategy as in \cite{JagNer84}, page~228,
where the case $a=1$ was proved.
Indeed, there it is shown that the intensity measure for individuals in
the $k$th
generation equals $\mu^{*k}$. Thus,
with $N^k$ denoting the set of individuals in generation $k$,
%
%e3.14 ###
\begin{eqnarray}
m_t^\chi(a)&=&\expec[z_t^\chi(a)]=\sum_{k=0}^{\infty
} a^k \expec\biggl[\sum_{i\in N^k}\chi(t-T_i) \biggr] \nonumber\\[-8pt]\\[-8pt]
&=&\int_0^{t} \chi(t-u)\sum_{k=0}^{\infty} a^k\mu^{*k}(du).\nonumber
\end{eqnarray}
For part (b), we follow
the identification of $\operatorname{Var}(z_t^\chi)$ in \cite{JagNer84}, Theorem 3.2 and
Corollary 3.3.
We use the covariance partition
%
%e3.15 ###
\begin{eqnarray}
\operatorname{Cov}(z_t^{\chi_1}(a_1), z_t^{\chi_2}(a_2))
&=&\operatorname{Cov} (\expec[z_t^{\chi_1}(a_1)\mid{\mathcal A}_0], \expec
[z_t^{\chi_2}(a_2)\mid{\mathcal A}_0] )\nonumber\\[-8pt]\\[-8pt]
&&{} +\expec[\operatorname{Cov}(z_t^{\chi_1}(a_1),z_t^{\chi_2}(a_2))\mid{\mathcal A}_0 ],\nonumber
\end{eqnarray}
where ${\mathcal A}_0$ is the $\sigma$-algebra generated by the lives
of the individuals in the first generation (the root is considered to
be in generation zero).\vadjust{\goodbreak} Then, the branching property of a CTBP gives that
%
%e3.16 ###
\begin{equation} z_t^{\chi}(a)=\chi(t)+a\sum_{j\colon G(j)=1}
z^{\chi}_{t-T_j}(j;a),
\end{equation}
where $ ((z^{\chi_1}_{t}(j;a_1),z^{\chi_2}_{t}(j;a_2)) )_{j, t\geq
0}$ is, conditionally on ${\mathcal A}_0$,\vspace*{1pt}
a sequence of i.i.d.\ random processes with law $((z^{\chi
_1}_{t}(a_1),z^{\chi_2}_{t}(a_2)))_{t\geq0}$. Therefore,
%
%e3.17 ###
\begin{equation}
\operatorname{Cov} (z_t^{\chi_1}(a_1),z_t^{\chi_2}(a_2)\mid
{\mathcal A}_0 ) =a_1a_2\sum_{j\dvtx G(j)=1} \mathrm{C}^{\chi_1,\chi
_2}_{t-T_j}(a_1,a_2),
\end{equation}
where we abbreviate
%
%e3.18 ###
\begin{equation} \mathrm{C}^{\chi_1,\chi_2}_{t}(a_1,a_2)=\operatorname{Cov}(z_t^{\chi_1}(a_1), z_t^{\chi_2}(a_2)).
\end{equation}
Thus,
%
%e3.19 ###
\begin{equation} \expec[\operatorname{Cov} (z_t^{\chi_1}(a_1),z_t^{\chi
_2}(a_2)\mid{\mathcal A}_0 ) ] =a_1a_2\int_0^t \mathrm{C}^{\chi_1,\chi
_2}_{t-v}(a_1,a_2)\mu(du).
\end{equation}
Further,
%
%e3.20 ###
\begin{equation} \expec[z_t^{\chi}(a)\mid{\mathcal A}_0]=\chi
(t)+a\int_0^t m^{\chi}_{t-u}(a) \PP(du),
\end{equation}
where $(\PP(t))_{t\geq0}$ is the intensity process of the first individual.
Therefore, we arrive at
%
%e3.21 ###
\begin{equation} \mathrm{C}^{\chi_1,\chi_2}_{t}(a_1,a_2)
=h_{a_1,a_2}(t)+(a_1a_2)\int_0^t \mathrm{C}^{\chi_1,\chi
_2}_{t-u}(a_1,a_2)\mu(du),
\end{equation}
where
%
%e3.22 ###
\begin{equation} h_{a_1,a_2}(t) = a_1a_2\operatorname{Cov} \biggl(\int_0^t
m_{t-u}^{\chi_1}(a_1) \PP(du), \int_0^t m_{t-u}^{\chi_2}(a_2)\PP
(du) \biggr).
\end{equation}
Iterating this equation yields (\ref{Cov-form}).

As before, for $\PP$ denoting the offspring distribution point process
[given by \eqref{eqn:pp-def}] and for every function $F\dvtx\R\to
\R$, note that
%
%e3.23 ###
\begin{equation} \label{cov-PPPa} \int_0^\infty F(x)\PP(dx)
\stackrel
{d}{=} f(\Pi),
\end{equation}
where $\Pi$ is a rate $1$ Poisson point process,
$f(x) \equiv F(x^s)$ and where the function $f$ applied to a point
process $\Pi$ is defined as $f(\Pi) \equiv\sum_{X\in\Pi} f(X)$.
By properties of functionals of the Poisson point process, we have that
%
%e3.24 ###
\begin{eqnarray}  \label{cov-PPPb}
\qquad \operatorname{Cov} \biggl(\int_0^t F_1(u)\PP(du),
\int_0^t F_2(u)\PP(du) \biggr) &=&\int_0^{t^{1/s}} F_1(u^s)F_2(u^s)\,du\nonumber\\[-8pt]\\[-8pt]
&=&s^{-1} \int_0^{t} u^{1/s-1} F_1(u)F_2(u)\,du.\nonumber
\end{eqnarray}
Therefore, we obtain
%
%e3.25 ###
\begin{equation} h_{a_1,a_2}(t)=\frac{a_1a_2}{s}\int_0^{t}
u^{1/s-1}m_{t-u}^{\chi_1}(a_1)m_{t-u}^{\chi_2}(a_2)\,du.\vspace*{-2pt}
\end{equation}
\upqed
\end{pf}

The following proposition (adapted mainly from \cite{JagNer84}, Theorem 3.5) captures
all we require to know about the asymptotics of the mean and variance
of a
single-vertex characteristic $\chi$.

%We recall that $\PP$ denoted the offspring point
%process. Define $\hat{\mu}(2\lambda)$ as the constant
% \eqn{
% \hat{\mu}(2\lambda) = \int_0^\infty\e^{-\lambda t}\mu(dt) < 1,
% }
%where the final inequality is due to the definition of
%$\lambda$ (replacing $2\lambda$ by $\lambda$ in the above
%equation gives us $1$). Define the random variable $\hat{\PP}(
% \eqn{
% \hat{\PP}(\lambda) = \int_0^\infty\e^{-\lambda x} \PP(dx).
% }
%By a computation as in \refeq{cov-PPPa}--\refeq{cov-PPPb}
%(with $t=\infty$),
% \eqn{
% \operatorname{Var}(\hat{\PP}(\lambda)) = \int_0^\infty\e^{-2\lambda x^s} dx <
% }
%This shall be crucial, as this constant appears in the limiting
%formula for the variance
%of the single-vertex characteristic below:

\begin{Proposition}[(Asymptotics of mean and variance for one-vertex
characteristics)]\label{prop:asymp-1-charac-exp}
For regular single vertex characteristics $\chi$, and all $a\geq0$,
\begin{longlist}[(a)]
\item[(a)] As $t\to\infty$,
%
%e3.26 ###
\begin{equation} \e^{-\lambda a^{s} t}\expec(z_t^\chi(a)) \to a^s
s\lambda\int_0^\infty\e^{-\lambda a^s y} \chi(y)\, dy.\vspace*{-2pt}
\end{equation}
When $a=a_t\rightarrow1$, then the convergence holds where in the
right-hand side
the value $a=1$ is substituted.

\item[(b)] As $t\to\infty$, when $a_1,a_2 \geq0$ with
$a^{s}_1+a_2^{s}-a_1^sa_2^s>0$,
%
%e3.27 ###
\begin{eqnarray}
&&\e^{-\lambda(a^{s}_1+a_2^{s}) t} \operatorname{Cov}(z_t^{\chi
_1}(a_1),z_t^{\chi_2}(a_2))\nonumber\\[-2pt]
&&\qquad  \to\frac{(a_1a_2)^{2+s}\lambda
^2s^2}{(a^{s}_1+a_2^{s})^{1/s}(a^{s}_1+a_2^{s}-a_1^sa_2^s)}\\[-2pt]
&&\qquad \quad
{}\times \int
_0^\infty\chi_1(x)\e^{-\lambda a_1^sx} \,dx \int_0^\infty\chi
_2(x)\e^{-\lambda a_2^sx} \,dx.\nonumber\vspace*{-2pt}
\end{eqnarray}
When $\vec{a}=\vec{a}_t\rightarrow(1,1)$, then the convergence holds
where in the right-hand sides
the value $\vec{a}=(1,1)$ is substituted.

\item[(c)] With $z_t=|\BP_t|$, there exists a random variable $W$ with $W>0$
a.s.\ such that $\e^{-\lambda t} z_t $ converges
almost surely and in $L^2$ to $W$ and further, for any single-vertex
regular characteristic,
%
%e3.28 ###
\begin{equation} \e^{-\lambda t} z_t^\chi\convas W \lambda\int
_0^\infty\chi(x)\e^{-\lambda x}\,dx,\vspace*{-2pt}
\end{equation}
and the convergence also holds in $L^2$.
\end{longlist}
\end{Proposition}
%
%carefully!}

\begin{pf}
Part (a) for $a\,{=}\,1$ is \cite{Jage75}, Theorem 5.2.8. For
$a\,{\neq}\,1$, we start~from
%
%e3.29 ###
\begin{equation} m_t^\chi(a) =\int_0^{t} \chi(t-u) \sum
_{k=0}^{\infty} a^k\mu^{*k}(du) =\int_0^t \chi(t-u) \tilde{\mu
}_{a}(du).\vspace*{-2pt}
\end{equation}
Define the measure with density $p_a(u)\,du$ via the equation $\e
^{-\lambda a^s u}\tilde{\mu}_{a}(du)=p_a(u)\,du +\e^{-\lambda a^s
u}\delta_0(du)$, then we obtain
%
%e3.30 ###
\begin{eqnarray}
\e^{-\lambda a^s t}m_t^\chi(a) &=&\e^{-\lambda a^s
t}\chi(t)+\int_0^t \chi(t-u)\e^{-\lambda a^s (t-u)} p_a(u)\,du \nonumber\\[-9pt]\\[-9pt]
&=&\e^{-\lambda a^s t}\chi(t)+\int_0^t \chi(v)\e^{-\lambda a^s v}
p_a(t-v)\,dv.\nonumber\vspace*{-2pt}
\end{eqnarray}
By Lemma \ref{lem-dens-tildemu}(a--b), we have that $p_a(u)$ is
uniformly bounded on $[1,\infty)$ and bounded by $c u^{1/s-1}$ on
$[0,1]$, while
and $p_a(u)\rightarrow a^s \lambda s$ when $u\rightarrow\infty$.
Thus, by dominated convergence,
%
%e3.31 ###
\begin{equation} \e^{-\lambda a^s t}m_t^\chi(a)\rightarrow a^s
\lambda s \int_0^\infty\chi(v)\e^{-\lambda a^s v}\,dv.
\end{equation}
The proof when $a_t\rightarrow1$ is identical.

See \cite{JagNer84}, Theorem 3.5, for parts (b) for $a_1=a_2=1$ and for
part (c).
For the proof of part (b) for $(a_1,a_2)\neq(1,1)$, we start with
(\ref{Cov-form})
and (\ref{ha-cov-def}). By part (a), we have that
%
%e3.32 ###
\begin{eqnarray}\label{ht-asy}
&&\e^{-\lambda(a^{s}_1+a_2^{s}) t}h_{a_1,a_2}(t)\nonumber \\
&&\qquad =\frac{a_1a_2}{s}\int_0^{t} u^{1/s-1} \e^{-\lambda(a^{s}_1+a_2^{s}) u}
\bigl(\e^{-\lambda a^{s}_1 (t-u)}m_{t-u}^{\chi_1}(a_1) \bigr) \nonumber\\[-8pt]\\[-8pt]
&&\qquad \quad\hphantom{\frac{a_1a_2}{s}\int_0^{t}}
{}\times\bigl(\e^{-\lambda
a_2^{s} (t-u)}m_{t-u}^{\chi_2}(a_2) \bigr)\,du\nonumber \\
&&\qquad \sim\frac
{a_1a_2}{s} M^{\chi_1}(a_1)M^{\chi_2}(a_2)\int_0^{\infty}
u^{1/s-1} \e^{-\lambda(a^{s}_1+a_2^{s})u}\,du,\nonumber
\end{eqnarray}
where we define
%
%e3.33 ###
\begin{equation} M^{\chi}(a)= a^s s\lambda\int_0^\infty\e
^{-\lambda y a^s} \chi(y) \,dy.
\end{equation}
Further, note that, by (\ref{lambda-s-comp}),
%
%e3.34 ###
\begin{eqnarray}
\int_0^{\infty} u^{1/s-1} \e^{-\lambda
(a^{s}_1+a_2^{s})u}\,du &=&(a^{s}_1+a_2^{s})^{-1/s} \lambda^{-1/s} \Gamma
(1/s)\nonumber\\[-8pt]\\[-8pt]
&=&s(a^{s}_1+a_2^{s})^{-1/s}.\nonumber
\end{eqnarray}
Then we rewrite
%
%e3.35 ###
\begin{eqnarray}
&&\e^{-\lambda(a^{s}_1+a_2^{s}) t} \operatorname{Cov}(z_t^{\chi_1}(a_1),z_t^{\chi_2}(a_2))\nonumber\\[-8pt]\\[-8pt]
&&\qquad  =\int_0^t \e^{-\lambda
(a^{s}_1+a_2^{s}) (t-u)} h_{a_1,a_2}(t-u) \e^{-\lambda
(a^{s}_1+a_2^{s}) u}\tilde{\mu}_{a_1a_2}(du).\nonumber
\end{eqnarray}
Now, by (\ref{p(u)-conv}), for $u$ large,
%
%e3.36 ###
\begin{equation} \e^{-\lambda(a^{s}_1+a_2^{s}) u}\tilde{\mu
}_{a_1a_2}(du) \sim(a_1a_2)^s \lambda s \e^{-\lambda
(a^{s}_1+a_2^{s}-a_1^sa_2^s) u}\,du,
\end{equation}
which is integrable, so that substitution of this asymptotics in
(\ref{Cov-form}) and using dominated convergence, proves that
%
%e3.37 ###
\begin{eqnarray}
\qquad &&\e^{-\lambda(a^{s}_1+a_2^{s}) t} \operatorname{Cov}\bigl(z_t^{\chi
_1}(a_1),z_t^{\chi_2}(a_2)\bigr) \nonumber \\
\qquad &&\qquad \sim (a_1a_2) M^{\chi
_1}(a_1)M^{\chi_2}(a_2)(a^{s}_1+a_2^{s})^{-1/s} \int_0^{\infty}
\e^{-\lambda(a^{s}_1+a_2^{s})
u}\tilde{\mu}_{a_1a_2}(du)\nonumber\\[-9pt]\\[-9pt]
\qquad&&\qquad\sim (a_1a_2)^2M^{\chi_1}(a_1)M^{\chi_2}(a_2)(a^{s}_1+a_2^{s})^{-1/s}
(a^{s}_1+a_2^{s}-a_1^sa_2^s)^{-1/s}\nonumber\\[-2pt]
\qquad &&\qquad =\frac
{(a_1a_2)^{1+s}s}{(a^{s}_1+a_2^{s})^{1/s}(a^{s}_1+a_2^{s}-a_1^sa_2^s)}
M^{\chi_1}(a_1)M^{\chi_2}(a_2).\nonumber
\end{eqnarray}
since, by Lemma \ref{lem-dens-tildemu},
\begin{eqnarray}
\int_0^{\infty} \e^{-\lambda b_1 u}\tilde{\mu}_{b_2}(du)
&=&
1+\int_0^{\infty} \e^{-\lambda (b_1-b_2^s)u}p_{b_2}(u)\,du\nonumber\\[-2pt]
&=&1+\int_0^{\infty} b_2^s \e^{-\lambda (b_1-b_2^s)u}p_{1}(u b_2^s)\,du\nonumber\\[-2pt]
&=&1+\int_0^{\infty}\e^{-\lambda
b_2^{-s}(b_1-b_2^s)u}p_{1}(u)\,du\nonumber\\[-10pt]\\[-10pt]
&=&\int_0^{\infty}\e^{-\lambda b_2^{-s}(b_1-b_2^s)u} \mu(du)
\nonumber\\[-2pt]
&=&\bigl(b_2^{-s}(b_1-b_2^s)\lambda\bigr)^{-1/s} \Gamma(1+1/s)\nonumber\\[-2pt]
&=&b_2(b_1-b_2^s)^{-1/s}\nonumber
\end{eqnarray}
by (\ref{lapl-trans}).
This proves\vspace*{1pt} the claim when $a^{s}_1+a_2^{s}-a_1^sa_2^s>0$.
When $\vec{a}=\vec{a}_t\rightarrow(1,1)$, then the above asymptotics
holds with
$\vec{a}=(1,1)$ substituted on the right-hand side\ since (\ref{ht-asy}) holds with
$\vec{a}=(1,1)$ substituted on its right-hand side.\vspace*{-3pt}
\end{pf}

%Let us start by stating two ramifications of this proposition:\\
%(I) First taking $\chi(\cdot) = 1$ and recalling that we used $z_t$ to
%denote the size of $\BP_t$ and using Lemma \ref{lemma:malthus}(d) on
%the form of $\mu^{*j}$ we get that first
% \eqn{
% \expec(z_t) = \int_0^t \sum_{j=0}^\infty\mu^{*j}(du) =\sum_{j=0}^
% }
%Using Proposition \ref{prop:asymp-1-charac-exp} we get the asymptotics
% \eqn{
% \e^{-\lambda t} \sum_{j=0}^\infty\frac{(\lambda t)^{\frac{j}{s}}}{
% \label{eqn:ident-e-sum}
% }
%(\mathit{II}) Next, taking $\chi(t) = t^{1/s-1}$ we get
% \eqn{
% \e^{-\lambda t} \int_0^\infty(t-u)^{1/s-1}\sum_{j=0}^\infty
% =\lambda s\frac{\Gamma(\frac{1}{s})}{\Gamma(\frac{1}{s}+1)}
% =\lambda s^2.\label{eqn:ident-e-sum-chi}}

%s3.3 ###
\subsection{\texorpdfstring{Almost sure convergence of two-vertex factor: Proof of
Theorem \protect\ref{theo:denom-as-convg}}
{Almost sure convergence of two-vertex factor: Proof of Theorem 2.5}}
\label{sec-as-two-vertex}
In this section, we prove Theorem \ref{theo:denom-as-convg}.
Throughout the proof, we shall abbreviate $p\,{=}\,1/s\,{-}\,1$.
Note that, for any fixed $0\,{<}\,\eps\,{<}\,B\,{<}\,\infty$, we can write~$\ztss$~as
%
%e3.38 ###
\begin{equation} \label{eqn:ztss-decomp}
\ztss=I^{(1)}_t(\eps,B)+ I^{(2)}_t(B)+ I^{(3)}_t(\eps),
\end{equation}
where
\begin{eqnarray*}
I^{(1)}_t(\eps,B) &=& \sum_{j\in\BP_t^{(2)}: \eps< t-T_j<B}
\sum_{i\in\BP_t^{(1)}} ([t-T_i]+ [t-T_j])^{p},\\[-2pt]
I^{(2)}_t(B) &=& \sum_{j\in\BP_t^{(2)}: t-T_j > B}
\sum_{i\in\BP_t^{(1)}} ([t-T_i]+ [t-T_j])^{p},\\[-2pt]
I^{(3)}_t(\eps) &= &\sum_{j\in\BP_t^{(2)}: t-T_j<\eps}
\sum_{i\in\BP_t^{(1)}} ([t-T_i]+ [t-T_j])^{p}.
\end{eqnarray*}
Thus to prove the result it is enough to show that for each fixed $\eps
,B$ we have
%
%e3.39 ###
\begin{eqnarray} \label{eqn:I1-conv}
&&\e^{-2\lambda t} I^{(1)}_t(\eps,B)\nonumber\\[-10pt]\\[-10pt]
&&\qquad \convas W^{(1)}
W^{(2)} \lambda^2\int_\eps^B\int_0^\infty(x_1+x_2)^{1/s-1} \e
^{-\lambda x_1} \e^{-\lambda x_2} \,dx_1\, dx_2,\nonumber
\end{eqnarray}\vspace*{-12pt}
\begin{eqnarray}
\label{eqn:I2-conv}
\limsup_{B\to\infty}\limsup_{t\to\infty}\e^{-2\lambda t} I^{ (2)}_t(B)&=& 0 \quad
\mbox{and}\nonumber\\[-8pt]\\[-8pt]
\limsup_{\eps\to0}\limsup_{t\to\infty}\e^{-2\lambda t} I^{
(3)}_t(\eps)&=&0.\nonumber
\end{eqnarray}
We shall start by proving \eqref{eqn:I1-conv}. The following lemma
shall be crucial in
our proof.
\begin{Lemma}[(Sup convergence of characteristics)]\label{lemma:sup-convg-char}
As $t\rightarrow\infty$,
%
%e3.41 ###
\begin{equation} \qquad \sup_{x\in[\eps,B]} \biggl|\e^{-\lambda t} \sum_{i\in
\BP_t^{(1)}} (x+ [t-T_i])^{p}- W^{(1)}\lambda\int_0^\infty(x+y)^{p}
\e^{-\lambda y} \,dy \biggr|\convas0,
\end{equation}
where $W^{(1)}$ is the almost sure limit of $\e^{-\lambda t}
z_t^{(1)}$.
\end{Lemma}

\begin{pf}
Consider the (random) functions
%
%e3.42 ###
\begin{equation} f_t(x) = \e^{-\lambda t} \sum_{i\in\BP_t^{(1)}}
(x+ [t-T_i])^{p},\qquad  x\in[\eps,B].
\end{equation}
Note that for $p<0$ these functions are monotonically decreasing, while
for $p>0$ they are increasing functions and they are all continuous
when defined on the compact interval
$[\eps,B]$. Further, for each fixed $x\in[\eps,B]$, by Proposition
\ref{prop:asymp-1-charac-exp}(c), pointwise we have on a set of
measure one,
%
%e3.43 ###
\begin{equation} f_t(x) \convas W^{(1)}\lambda\int_0^\infty(x+y)^{p}
\e^{-\lambda y} \,dy.
\end{equation}
Thus, to show the a.s.\ sup convergence, by the Arzela--Ascoli theorem
(see, e.g.,\ \cite{rudin-real}),
it is enough to show that the above family of functions are a.s.\
equicontinuous, that is, for any $x\in[\eps,B]$ and any given $\delta>
0$ there exists $\eta(x)> 0$ independent of $t$ such that for all $t$:
%
%e3.44 ###
\begin{equation}\label{eqn:arz-asc}
|f_t(x) -f_t(y)| < \delta\qquad \mbox{for all } y\in
[x-\eta(x),x+\eta(x)]\cap[\eps,B].
\end{equation}
We separate between the cases $p<1$ and $p\geq1$.

\textit{Case} 1: $p< 1$. In this case note that for any $l_1, l_2>
0$ and $a> 0$, we have
\begin{eqnarray} \label{eqn:pleq1}
|(l_1+a)^p - (l_2+a)^p| &=&|p-1|\int_{l_1}^{l_2} (x+a)^{p-1}\,dx\nonumber\\[-8pt]\\[-8pt]
&\leq&|p-1|\int_{l_1}^{l_2} x^{p-1}\,dx \qquad \mbox{since } p-1< 0.\nonumber
\end{eqnarray}
By the continuity of the function $g(x) = x^{p}$, for any $x\in[\eps
,B]$ and $\delta^\prime>0$, we can choose $\eta^\prime(x)$ small
such that for $y\in[\eps,B], |y-x|< \eta^\prime(x)$ we have
%
%e3.45 ###
\begin{equation}
|x^p -y^p |< \delta^\prime.
\end{equation}
This implies from \eqref{eqn:pleq1} applied individually to the
functions $g_i(x) = (x+[t-T_i])^{1/s-1}$ for $i\in\BP_t^{(1)}$ that
\[
|f_t(x) -f_t(y)| < \delta^\prime\e^{-\lambda t} z_t^{(1)},
\]
where we recall that $z_t^{(1)} = |\BP^{(1)}_t|$. Since $\e
^{-\lambda t} z_t^{(1)}$ converges a.s.\ and\break $(\e^{-\lambda t}
z_t^{(1)})_{t\geq0}$ is bounded a.s., we obtain that,
on a set $A$ of measure one, for each $\omega\in A$, we can find a
$\kappa(\omega)$ depending on the sample point $\omega$ but
independent of $t$, such that
\[
\sup_{t} \e^{-\lambda t} z_t^{(1)}(\omega) < \kappa(\omega).
\]
Now choosing $\delta^{\prime} = \delta/\kappa(\omega)$ gives us a
$\eta(x) = \eta(x,\omega)$ such that
\eqref{eqn:arz-asc} is satisfied. This proves the result for $p<1$.

\textit{Case} 2: $p\geq1$. Here note that for any $a>0$ and
$x,y\in[\eps,B]$,
we have, by the mean value theorem
\[
|(x+a)^{p}-(y+a)^p| \cases{= p(z+ a)^{p-1} |y-x|,&\quad  $a\in[x,y]$,\vspace*{2pt}
\cr
\leq p(B+a)^{p-1} |y-x|,&\quad  since $p-1\geq0$.}
\]
This implies that, for $x,y \in[\eps,B]$,
%
%e3.46 ###
\begin{equation} |f_t(x) -f_t(y)| < H_t |x-y|,
\end{equation}
where, by Proposition \ref{prop:asymp-1-charac-exp}(c),
%
%e3.47 ###
\begin{eqnarray}
H_t &=& \e^{-\lambda t} p\sum_{i\in\BP^{(1)}_t} (B+
[t-T_i])^{p-1} \nonumber\\[-8pt]\\[-8pt]
&\convas& p W^{(1)}\int_0^\infty(B+y)^{p-1} \lambda\e
^{-\lambda y} \,dy.\nonumber
\end{eqnarray}
This proves that \eqref{eqn:arz-asc} holds also when $p\geq1$,
and completes the proof of Lemma~\ref{lemma:sup-convg-char}.
\end{pf}

\begin{pf*}{Completion of the proof of (\ref{eqn:I1-conv})}
Write
%
%e3.48 ###
\begin{equation}
\BP_t^{(2)}(\eps,B)= \bigl\{v\in\BP^{(2)}_t\dvtx\eps
< t-T_j<B \bigr\}.
\end{equation}
Then we have
%
%e3.49 ###
\begin{eqnarray}\label{eqn:diff}
\qquad &&\biggl|I_t^{(1)}(\eps,B)- \e^{-\lambda t}\sum_{j\in\BP_t^{
(2)}(\eps,B) }W^{(1)}\int_0^\infty
([t-T_j]+y)^{1/s-1} \lambda\e^{-\lambda y} \,dy \biggr| \nonumber\\[-8pt]\\[-8pt]
\qquad &&\qquad \leq Q_1(t)\e^{-\lambda t} z_t^{(2)},\nonumber
\end{eqnarray}
where
%
%e3.50 ###
\begin{eqnarray}
Q_1(t) &=& \sup_{x\in[\eps,B]} \biggl|\e^{-\lambda t}\sum
_{i\in\BP_t^{(1)}}(x+[t-T_i])^{1/s-1} \nonumber\\[-8pt]\\[-8pt]
&&\hphantom{\sup_{x\in[\eps,B]} \biggl|}
{}- W^{(1)}\int_0^\infty
(x+y)^{1/s-1} \lambda\e^{-\lambda y} \,dy \biggr|.\nonumber
\end{eqnarray}
Lemma \ref{lemma:sup-convg-char} now implies that the term on the
right-hand side\ of \eqref{eqn:diff} converges to 0 a.s. Thus, to complete the
proof, it is enough to show that
\begin{eqnarray*}
&&\e^{-\lambda t}\sum_{j\in\BP_t^{(2)}(\eps,B) }W^{
(1)}\int_0^\infty([t-T_j]+y)^{p} \lambda\e^{-\lambda y} \,dy \\
&&\qquad \convas
W^{(1)} W^{(2)} \lambda^2\int_\eps^B\int_0^\infty(x_1
+x_2)^{p} \e^{-\lambda(x_1+x_2)}\,dx_1 \,dx_2.
\end{eqnarray*}
This follows by taking the characteristic
%
%e3.51 ###
\begin{equation} \chi_2(a)=
\cases{
\displaystyle \int_0^\infty(a+x_2)^{p} \e^{-\lambda x_2} \,dx_2, &\quad
if $\eps\leq a\leq B$, \vspace*{2pt}\cr
0, &\quad  if $a\notin[\eps,B]$
}
\end{equation}
and using Proposition \ref{prop:asymp-1-charac-exp}(c) for the
branching process $\BP_t^{(2)}$.
\end{pf*}

\begin{pf*}{Completion of the proof of (\ref{eqn:I2-conv})}
First, consider the term $I_t^{(3)}(\eps)$. Note that
%
%e3.52 ###
\begin{equation} I_t^{(3)}(\eps) \leq\bigl[z_t^{(2)}-z_{t-\vep}^{ (2)}\bigr]
\sum_{j\in\BP_t^{(1)}} (\eps+ t-T_j)^{p}.
\end{equation}
By Proposition \ref{prop:asymp-1-charac-exp}(c)
%
%e3.53 ###
\begin{eqnarray}
&&\e^{-2\lambda t} \bigl[z_t^{(2)}-z_{t-\vep}^{(2)}\bigr]\sum
_{j\in\BP_t^{(1)}} (\eps+ t-T_j)^{p} \nonumber\\[-8pt]\\[-8pt]
&&\qquad \convas\bigl[W^{(2)} [1-\e
^{-\lambda\vep}] \bigr]\cdot\biggl[W^{(1)} \int_0^\infty(\eps+y)^{p} \e
^{-\lambda y} \,dy \biggr]\convas0,\nonumber
\end{eqnarray}
when $\vep\downarrow0.$
%Thus
% \eqn{
% \lim_{\eps\to0} [W^{(2)}\int_0^\eps\e^{-\lambda x} dx ]\cdot
%[W^{(1)} \int_0^\infty(\eps+y)^{p} \e^{-\lambda y} dy ]\convas0.
% }
This proves the last convergence result in \eqref{eqn:I2-conv}.

To prove the first convergence result in \eqref{eqn:I2-conv},
note that arguing as in the proof of \eqref{eqn:I1-conv}, we have for
all $p$ and $x_1, x_2> 0$,
%
%e3.54 ###
\begin{equation} (x_1+x_2)^{p} \leq(2^{p}\vee1)(x_1^{p}+ x_2^{p}),
\end{equation}
where $a\vee b=\max\{a,b\}$. Thus,
%
%e3.55 ###
\begin{equation} I_t^{(2)}(B) \leq(2^{p}\vee1)\bigl[z_t^{
{(1)}}z_{t-B}^{\chi_{B},(2)}+z_t^{\chi,{(1)}}z_{t-B}^{(2)}\bigr],
\end{equation}
where $\chi_B(x)=(B+x)^p, \chi(x)=x^p$. Now again, by Proposition
\ref{prop:asymp-1-charac-exp}(c),
%
%e3.56 ###
\begin{eqnarray}
&&\e^{-\lambda t} \bigl[z_t^{{(1)}}z_{t-B}^{\chi_{B},
(2)}+z_t^{\chi,{(1)}}z_{t-B}^{(2)}\bigr] \nonumber\\
&&\qquad \convas W^{(1)} W^{(2)} \e
^{-\lambda B} \biggl[\int_0^\infty\int_0^\infty\e^{-\lambda x} (B+x)^{p} \e^{-\lambda
x} \,dx\,dy\\
&&\qquad  \hspace*{111pt}
{}+ \int_0^\infty\int_0^\infty\e^{-\lambda x} x^{p} \e^{-\lambda x} \,dx \,dy\biggr],\nonumber
\end{eqnarray}
which converges a.s.\ to 0 when $B\rightarrow\infty$.
This completes the proof of \eqref{eqn:I2-conv}.
\end{pf*}

%s3.4 ###
\subsection{Mean and variance of two-vertex characteristic}
\label{sec-two-vertex}

In this section, we shall analyze two-vertex characteristics. This sets
the stage for the
proof of the asymptotics for the hopcount in Theorem \ref
{theo:two-vertex-charac-ratio}. Define, for $\vec{a}=(a_1,a_2)$,
%
%e3.57 ###
\begin{equation} \label{z(1,2)-def} z_t^{(1,2)}(\vec{a}) = \sum
_{i\in
\BP^{(1)}_t} \sum_{j\in\BP^{(2)}_t}
a_1^{G^{(1)}(i)}a_2^{G^{(2)}(j)}([t-T_i]+ [t-T_j])^{p},
\end{equation}
where we recall that $G^{(i)}(v)$ is the \emph{generation} of
vertex $v\in\BP^{(i)}_t$.

\begin{Lemma}[(Expectation and variance of two-vertex characteristics)]\label{lem-mean-var-two-vertex-charac}
Consider two independent CTBPs
$\BP^{(1)}_t \mbox{ and } \BP_t^{(2)}$ driven by the
offspring distribution $\PP$. Then
\begin{longlist}[(a)]
\item[(a)]
%
%e3.58 ###
\begin{equation}
\expec\bigl[z_t^{(1,2)}(\vec{a})\bigr]= \int_0^t\int_0^t
([t-v]+[t-u])^{p}\tilde{\mu}_{a_1}(dv) \tilde{\mu}_{a_2}(du).
\end{equation}
\item[(b)]
%
%e3.59 ###
\begin{eqnarray}
&&\operatorname{Cov}\bigl(z_t^{(1,2)}(\vec{a}), z_t^{(1,2)}(\vec
{b})\bigr)\nonumber\\[-8pt]\\[-8pt]
&&\qquad =\int_0^t h^{(1)}_{\vec{a},\vec{b}}(t-u,t) \tilde{\mu
}_{a_2b_2}(du)
 +\int_0^t h^{(2)}_{\vec{a},\vec{b}}(t-u,t)\tilde{\mu}_{a_1b_1}(du),\nonumber
\end{eqnarray}
where
%
%e3.61 ###
%e3.60 ###
\begin{eqnarray}
h_{\vec{a},\vec{b}}^{(1)}(v,t) &=&\frac{a_2b_2}{s}\int
_0^{v} \int_0^{v-u} \int_0^{v-u} \int_0^t \int_0^t u^{p}
([t-u_2]+[v-u-u_1])^{p}\nonumber\\
&&\hspace*{123pt}{} \times
([t-v_2]+[v-u-v_1])^{p}\nonumber\\[-8pt]\\[-8pt]
&&\hspace*{123pt}{} \times
\tildmu_{a_1}(du_1)\tildmu_{b_1}(dv_1)\tildmu_{a_2}(du_2)\nonumber \\
&&\hspace*{123pt}{} \times\tildmu
_{b_2}(dv_2)\,du,\nonumber\\
\qquad \quad h^{(2)}_{\vec{a},\vec{b}}(v,t)&=&\frac
{a_1b_1}{s}\int_0^{v} u^{p}\expec\bigl[z_{v-u}^{\tilde\chi^{(2)}_{t,
v-u,a_1}}(a_2)z_{v-u}^{\tilde{\chi}^{(2)}_{t,v-u,b_1}}(b_2) \bigr]\,du
\end{eqnarray}
with $p=1/s-1$ and
%
%e3.62 ###
%e3.76 ###
\begin{equation} \tilde{\chi}_{t,r,a_2}^{(2)}(x)=\int_0^r
(t-u_2+x)^{p}\tilde{\mu}_{a_2}(du_2).
\end{equation}
\end{longlist}
\end{Lemma}

\begin{pf}
%For part (a), we recall that for a given single-vertex
%characteristic $\chi$, $z_{t}^{\chi, (i)}$ was used to denote:
% \eqn{
% z_{t}^{\chi, {(i)}}(a) = \sum_{v\in\BP_t^{(i)}} a^{G^{
%(i)}(v)}\chi(t-T_v).
% }
%For the rest of the treatment, for ease of exposition we shall write
%the measure
% \eqn{
% \tildmu_a(du) = \sum_{j=0}^\infty a^j\mu^{*j}(du).
% }
We shall prove part (a) by conditioning on $\BP_t^{(1)}$.
Note we can write $z_t^{(1,2)}(\vec{a})=z_{t}^{\chi
_{t,a_1}^{(1)}, {(2)}}(a_2),$ where
%
%e3.63 ###
\begin{equation} \label{chi-1-a1} \chi_{t,a_1}^{(1)}(x)=\sum_{j\in
\BP^{(1)}_t} a_1^{G^{ (1)}(j)}(x+ [t-T_j])^{p}.
\end{equation}
Conditionally on $\BP_t^{(1)}$, the characteristic $\chi
_{t,a_1}^{(1)}$ is
deterministic. Therefore, Proposition \ref{prop:form-1-charac-exp}(a)
implies that
%
%e3.64 ###
\begin{eqnarray}  \label{cond-expec-z12}
\expec\bigl(z_t^{(1,2)}(\vec
{a})\mid\BP_t^{(1)}\bigr) &=& \sum_{j\in\BP_t^{(1)}} a_1^{G^{(1)}(j)} \int
_0^t ([t-u]+[t-T_j])^{p}\tildmu_{a_2}(du) \nonumber\\[-8pt]\\[-8pt]
&=& z_{t}^{\chi^{(2)}_{t,a_2}, {(1)}}(a_1),\nonumber
\end{eqnarray}
where $\chi^{(2)}_{t,a_2}$ is the characteristic
%
%e3.65 ###
\begin{equation} \label{chi-2-a2}
\chi^{(2)}_{t,a_2}(v)= \int_0^t
([t-u]+v)^{p}\tildmu_{a_2}(du).
\end{equation}
We complete the proof by noting that, for all $r$,
%
%e3.66 ###
\begin{eqnarray}  \label{m-chi-(2)-ident} m_{r}^{\chi
^{(2)}_{t,a_2}}(a_1) &=&\expec\bigl(z_{r}^{\chi^{(2)}_{t,a_2},
{(2)}}(a_1)\bigr)\nonumber\\[-8pt]\\[-8pt]
&=& \int_0^r\int_0^t ([r-v]+[t-u])^{p}\tildmu_{a_1}(dv)\tildmu
_{a_2}(du).\nonumber
\end{eqnarray}
Taking $r=t$ proves the claim in part (a).

For part (b), we again condition on $\BP_t^{(1)}$, for which we
use the covariance partition
%
%e3.67 ###
\begin{eqnarray}\label{eqn:var-decomp}
&&\operatorname{Cov}\bigl(z_t^{(1,2)}(\vec{a}), z_t^{(1,2)}(\vec
{b})\bigr)\nonumber \\
&&\qquad = \operatorname{Cov} \bigl(\expec\bigl(z_t^{(1,2)}(\vec{a})\mid\BP_t^{(1)}\bigr),
\expec\bigl(z_t^{(1,2)}(\vec{b})\mid\BP_t^{(1)}\bigr) \bigr) \nonumber\\[-8pt]\\[-8pt]
&&\qquad \quad \hspace*{35pt}
{}+ \expec\bigl(\operatorname{Cov}\bigl(z_t^{(1,2)}(\vec{a}), z_t^{(1,2)}(\vec{b})\mid\BP_t^{(1)}\bigr)
\bigr)\nonumber \\
&&\qquad =(I)_t+(\mathit{II})_t. \nonumber
\end{eqnarray}
Let us now tackle each of these two terms separately.

\textit{Term} $(I)_t$: For $(I)_t$, we use
the explicit formula for $\expec(z_t^{(1,2)}(\vec{a})\mid\BP
_t^{(1)})$ in~(\ref{cond-expec-z12}) and $\chi^{(2)}_{t,a_2}$ in (\ref{chi-2-a2})
to obtain that
%
%e3.68 ###
\begin{equation} (I)_t= \operatorname{Cov} \bigl(z_t^{\chi
^{(2)}_{t,a_2},{(1)}}(a_1), z_t^{\chi^{(2)}_{t,b_2}, {(1)}}(b_1) \bigr).
\end{equation}
Now using Proposition \ref{prop:form-1-charac-exp}(b), we get
%
%e3.69 ###
\begin{equation}\label{eqn:var-convolv}
(I)_t= \int_0^t \int_0^t h_{\vec{a},\vec
{b}}^{(1)}(t-u,t) \tilde{\mu}_{a_2b_2}(du),
\end{equation}
where $v\mapsto h_{\vec{a},\vec{b}}^{(1)}(v,t)$ is the function
%
%e3.70 ###
\begin{eqnarray} \label{hvecavecb-def}
h_{\vec{a},\vec{b}}^{(1)}(v,t)
&=&\frac{a_2b_2}{s}\int_0^{v} u^{p}m_{v-u}^{\chi^{
(2)}_{t,a_2}}(a_1)m_{v-u}^{\chi^{(2)}_{t,b_2}}(b_1)\,du\nonumber\\
&=&\frac
{a_2b_2}{s}\int_0^{v} u^{p} \int_0^{v-u} \int_0^{v-u}\chi
^{(2)}_{t,a_2}(v-u-u_1)\chi^{(2)}_{t,b_2}(v-u-v_1)\nonumber \\
&&\hspace*{106pt}{}\times \tildmu
_{a_1}(du_1)\tildmu_{b_1}(dv_1)\,du\nonumber\\[-8pt]\\[-8pt]
&=&\frac{a_2b_2}{s}\int
_0^{v} \int_0^{v-u} \int_0^{v-u} \int_0^t \int_0^t u^{p}
([t-u_2]+[v-u-u_1])^{p}\nonumber\\
&&\hspace*{123pt}
{}\times([t-v_2]+[v-u-v_1])^{p}
\tildmu_{a_1}(du_1)\nonumber\\
&&\hspace*{123pt}
{} \times\tildmu_{b_1}(dv_1)\tildmu_{a_2}(du_2)\tildmu
_{b_2}(dv_2)\,du.\nonumber
\end{eqnarray}

\textit{Term} $(\mathit{II})_t$: We again use that, conditionally on $\BP
_t^{(1)}$,
$z_t^{(1,2)}(\vec{a})=\break z_{t}^{\chi_{t,a_1}^{(1)}, {
(2)}}(a_2),$ where
$\chi_{t,a_1}^{(1)}$ was defined in (\ref{chi-1-a1}).
Therefore, we can again use Proposition~\ref
{prop:form-1-charac-exp}(b) to write
%
%e3.71 ###
\begin{equation} \label{eqn:htu-random}
\operatorname{Cov} \bigl(z_t^{(1,2)}(\vec{a}), z_t^{(1,2)}(\vec
{b})\mid\BP_t^{(1)}\bigr) =\int_0^t g_{\vec{a},\vec{b}}(t-u,t)\tildmu
_{a_2b_2}(du),
\end{equation}
where
%
%e3.72 ###
\begin{equation}
g_{\vec{a},\vec{b}}(v,t) =\frac{a_2b_2}{s}\int_0^{v}
u^{p}m_{v-u}^{\chi_{t,a_1}^{ (1)}}(a_2)m_{v-u}^{\chi
_{t,a_1}^{(1)}}(b_1)(b_2)\,du.
\end{equation}
Therefore,
%
%e3.73 ###
\begin{equation} (\mathit{II})_t= \int_0^t h^{(2)}_{\vec{a},\vec
{b}}(t-u,t)\tildmu_{a_2b_2}(du),
\end{equation}
where
%
%e3.74 ###
\begin{equation} h^{(2)}_{\vec{a},\vec{b}}(v,t)= \frac
{a_2b_2}{s}\int_0^{v} u^{p}\expec\bigl[m_{v-u}^{\chi_{t,a_1}^{
(1)}}(a_2)m_{v-u}^{\chi_{t,b_1}^{(1)}}(b_2)\bigr]\,du.
\end{equation}
We can now rewrite
%
%e3.75 ###
\begin{eqnarray}
m_{r}^{\chi_{t,a_1}^{(1)}}(a_2) &=&\int_0^{r} \sum
_{j\in\BP_t^{(1)}}a_1^{G^{ (1)}(j)}(t-u_2+t-T_j)^{p} \tilde{\mu
}_{a_2}(du_2)\nonumber \\
&=&\sum_{j\in\BP_t^{(1)}}a_1^{G^{(1)}(j)} \int_0^{r}
(t-u_2+t-T_j)^{p}\tilde{\mu}_{a_2}(du_2)\\
&=&z_t^{\tilde\chi
_{t,r,a_2}^{(2)}, {(1)}}(a_1),\nonumber
\end{eqnarray}
where
\begin{equation} \tilde{\chi}_{t,r,a_2}^{(2)}(x)=\int_0^r
(t-u_2+x)^{p}\tilde{\mu}_{a_2}(du_2).
\end{equation}
This completes the proof.
\end{pf}

%Note that in the above expressions, the variance is taken only over
%the disorder in $\PP$, namely in the context of the above, $z_{l,2}^{
% \[z_{r,2}^{\chi_2}= \sum_{j\in\BP^{(2)}_t} \int_0^r (r-u+
%(t-T_j))^{1/s-1}\tildmu(du).\]
%
%Now for any fixed $u,v$ we have that
% \[z_{t-u-v,2}^{\chi_2}= \int_0^{t-u-v} \sum_{j\in\BP^{(2)}_t}
%(t-u-v-x+ (t-T_j))^{1/s-1}\tildmu(dx) \]
%A truncation argument and Proposition \ref{prop:asymp-1-charac-exp}(c)
%implies that
% \[e^{-(2t-u-v)} z_{t-u-v,2}^{\chi_2} \stackrel{a.s., \mathbb{L}^2}{
%In particular this implies that for each fixed $u$
% \[e^{-4\lambda(t-u)} h(t-u) \stackrel{a.s., \mathbb{L}^2}{
%Using this and taking expectations in \eqref{eqn:htu-random} shows that
% \[e^{-4\lambda t} \expec(\Var(\ztss|\BP_t^{(2)}))\to[\Psi_1]^2
% \]
%This completes the proof.

\begin{Lemma}[(Asymptotics of mean and variance of two-vertex characteristics)]\label{lemma:mean-two-vertex-charac}
Consider two independent CTBPs $\BP^{(1)}_t $ and $\BP
^2_t$ driven by the offspring distribution $\PP$. Then
\begin{longlist}[(a)]
\item[(a)]
%
%e3.77 ###
\begin{eqnarray}  \label{asy-mean-two-vertex}
&&\e^{-\lambda(a_1^s+a_2^s) t}\expec\bigl(z_t^{(1,2)}(\vec{a})\bigr)\nonumber\\[-8pt]\\[-8pt]
&&\qquad  \to(\lambda s)^2 \int
_0^\infty\int_0^\infty(x_1+x_2)^{p}\e^{-\lambda(a_1^s x_1+a_2^s
x_2)}\,dx_1 \,dx_2.\nonumber
\end{eqnarray}
When $\vec{a}=\vec{a}_t\rightarrow(1,1)$, then the convergence holds
where in the right-hand sides
the value $\vec{a}=(1,1)$ is substituted.

\item[(b)] When $a^{s}_1+a_2^{s}-a_1^sa_2^s>0$ and
$b^{s}_1+b_2^{s}-b_1^sb_2^s>0$, there exists a constant~$A_{\operatorname{Cov}}(\vec{a}, \vec{b})$ such that
%
%e3.78 ###
\begin{equation}
\e^{-\lambda[(a_1^s+a_2^s)+(b_1^s+b_2^s)] t}\operatorname{Cov}\bigl(z_t^{ (1,2)}(\vec{a}), z_t^{(1,2)}(\vec{b})\bigr) \to A_{\operatorname{Cov}}(\vec{a}, \vec{b}).
\end{equation}
When $\vec{a}=\vec{a}_t\rightarrow(1,1), \vec{b}=\vec
{b}_t\rightarrow(1,1)$, then the convergence holds
where the right-hand side\ is replaced with $A_{\operatorname{Cov}}(\vec{1}, \vec{1})$.
\end{longlist}
\end{Lemma}

\begin{pf}
By Lemma \ref{lem-mean-var-two-vertex-charac}(a),
%
%e3.79 ###
\begin{eqnarray}
\quad &&\e^{-\lambda(a_1^s+a_2^s) t}\expec\bigl(z_t^{(1,2)}(\vec
{a})\bigr)\nonumber  \\
\quad &&\qquad =\e^{-\lambda(a_1^s+a_2^s) t}\int_0^t\int_0^t
([t-u]+[t-v])^{1/s-1} \tilde{\mu}_{a_1}(du)\tilde{\mu}_{a_2}(dv)\nonumber \\
\quad &&\qquad =\int_0^t\int_0^t \e^{-\lambda
[a_1^s(t-u)+a_2^s(t-v)]}([t-u]+[t-v])^{1/s-1} p_{a_1}(u) p_{a_2}(v)\,du
\,dv\\
\quad &&\qquad\quad {} +o(1)\nonumber\\
\quad &&\qquad =\int_0^t\int_0^t \e^{-\lambda[a_1^s
u+a_2^sv]}(u+v)^{1/s-1} p_{a_1}(t-u) p_{a_2}(t-v)\,du \,dv+o(1),\nonumber
\end{eqnarray}
where the $o(1)$ originates from part $\delta_{u,0}(du)$ in the
decomposition\break
$\e^{-\lambda a^s u}\tilde{\mu}_{a}(du)=p_a(u)\,du +\e^{-\lambda a^s
u}\delta_0(du)$,
and $\delta_{u,0}$ is the Dirac measure at $u=0$.
Now we again use Lemma \ref{lem-dens-tildemu}(a) and (b) together with
dominated
convergence to obtain the claim in part (a). The proof of part (b) is
similar, and we omit the details.
\end{pf}

\subsection{\texorpdfstring{CLT for two-vertex characteristic: Proof of Theorem \protect\ref{theo:two-vertex-charac-ratio}}
{CLT for two-vertex characteristic: Proof of Theorem 2.7}}\label{sec-pf-CLT}

In this section, we use and extend the theory developed in the previous
section to prove Theorem~\ref{theo:two-vertex-charac-ratio}. The plan
is as follows: We shall start by proving the result
for the CTBP in the summation instead of $\SWG_t^{(i)}$ and then
argue that the difference is negligible. The result is formulated as follows.

\begin{Theorem}[(CLT for two-vertex characteristic for CTBP)]\label{theo:two-vertex-charac-ratio-CTBP}
The two-vertex characteristic satisfies that, as $t\rightarrow\infty$,
%
%e3.80 ###
\begin{equation} \label{conv-two-vertex-CTBP} \frac{\sum_{i\in\BP
_t^{(1)}(x)} \sum_{j\in\BP_t^{ (2)}(y)}\chi_{ij}(t)}{\sum_{i\in
\BP_t^{(1)}} \sum_{j\in\BP_t^{(2)}}\chi_{ij}(t)} \convp\Phi
(x)\Phi(y).
\end{equation}
\end{Theorem}

%The remainder of this section is organized as follows. We start by
%proving Theorem \ref{theo:two-vertex-charac-ratio-CTBP},
%and after the completion of the proof, we shall complete the
%proof of Theorem \ref{theo:two-vertex-charac-ratio}, which shall
%be relatively easy, as there is little distinction on the
%appropriate time scale between the SWG and the CTBP.

\begin{pf}
Theorem \ref{theo:two-vertex-charac-ratio-CTBP} follows when we
show that,
writing $\vec{a}_t=(\e^{\alpha_1/\sqrt{s^2\lambda t}},\allowbreak\e^{\alpha
_2/\sqrt{s^2\lambda t}})$,
for some vector $\vec{\alpha}=(\alpha_1,\alpha_2)$,
%
%e3.81 ###
\begin{equation} \label{conv-two-vertex-CTBP-rep} \frac
{z_t^{(1,2)}(\vec{a}_t)\e^{-(\alpha_1+\alpha_2)\sqrt{\lambda t}}}
{z_t^{(1,2)}} \convp\e^{\alpha_1^2/2+\alpha_2^2/2}.
\end{equation}
Indeed, define the (random) measure $P$ on pairs $(X,Y)$ by
%
%e3.82 ###
\begin{equation} P(X\leq x, Y\leq y) =\frac{\sum_{i\in\BP
_t^{(1)}(x)} \sum_{j\in\BP_t^{ (2)}(y)}\chi_{ij}(t)}{\sum_{i\in
\BP_t^{(1)}} \sum_{j\in\BP_t^{(2)}}\chi_{ij}(t)}.
\end{equation}
Then, (\ref{conv-two-vertex-CTBP}) states that the pair $(X,Y)$ converges
in distribution to a pair of independent standard normal distributions,
which, in turn, follows when, for each $(\alpha_1,\alpha_2)\in\R$,
we have that
%
%e3.83 ###
\begin{equation} E[\e^{\alpha_1X+\alpha_2Y}]= \frac
{z_t^{(1,2)}(\vec{a}_t) \e^{-(\alpha_1+\alpha_2)\sqrt{\lambda
t}}}{z_t^{(1,2)}}\convp\e^{\alpha_1^2/2+\alpha_2^2/2},
\end{equation}
where $E$ denotes the expectation w.r.t.\ $P$.

In order to show (\ref{conv-two-vertex-CTBP-rep}), we show that
%
%e3.84 ###
\begin{equation} \label{aim-CLT-two-vertex} \e^{-2\lambda t}
\bigl(z_t^{(1,2)}(\vec{a}_t)\e^{-(\alpha_1+\alpha_2)\sqrt{\lambda t}}
-\e^{\alpha_1^2/2+\alpha_2^2/2}z_t^{(1,2)} \bigr)\convp0.
\end{equation}
Together with Theorem \ref{theo:denom-as-convg}, this then implies the
result, as
$\e^{-2\lambda t}z_t^{(1,2)}$ converges a.s.\ to a strictly
positive random variable.
Thus, we are left to pro\-ve~(\ref{aim-CLT-two-vertex}). We shall show
that the convergence in (\ref{aim-CLT-two-vertex}) in fact holds in~$L^2$.
For this, it is immediate that it suffices to study
%
%e3.86 ###
%e3.85 ###
\begin{eqnarray}
M_t(\vec{\alpha})&\equiv&\expec\bigl[z_t^{(1,2)}(\vec
{a}_t)\e^{-(\alpha_1+\alpha_2)\sqrt{\lambda t}} -\e^{\alpha
_1^2/2+\alpha_2^2/2}z_t^{(1,2)} \bigr],\\
Q_t(\vec{\alpha})&\equiv&\operatorname{Var}\bigl(z_t^{(1,2)}(\vec{a}_t)^2\bigr),\qquad
C_t(\vec{\alpha})\equiv\operatorname{Cov}\bigl(z_t^{(1,2)},z_t^{ (1,2)}(\vec{a}_t)\bigr).
\end{eqnarray}
In terms of these quantities, we can rewrite
%
%e3.87 ###
\begin{eqnarray}
&&\expec\bigl[ \bigl(z_t^{(1,2)}(\vec{a}_t)\e^{-(\alpha_1+\alpha
_2)\sqrt{\lambda t}} -\e^{\alpha_1^2/2+\alpha_2^2/2}z_t^{(1,2)}
\bigr)^2 \bigr]\nonumber\\
&&\qquad  =M_t(\vec{\alpha})^2+Q_t(\vec{\alpha})\e^{-2(\alpha
_1+\alpha_2)\sqrt{\lambda t}} +\e^{\alpha_1^2+\alpha_2^2}Q_t(\vec
{0})\\
&&
\qquad \quad {}- 2C_t(\vec{\alpha})\e^{-(\alpha_1+\alpha_2)\sqrt{\lambda t}}
\e^{\alpha_1^2/2+\alpha_2^2/2}.\nonumber
\end{eqnarray}
Therefore, we shall prove that
%
%e3.88 ###
\begin{equation}\label{aim-1}
\e^{-2\lambda t}M_t(\vec{\alpha})=o(1)
\end{equation}
and
%
%e3.89 ###
\begin{eqnarray} \label{aim-2}
&&\e^{-4\lambda t} \bigl( Q_t(\vec{\alpha})\e^{-2(\alpha
_1+\alpha_2)\sqrt{\lambda t}} +\e^{\alpha_1^2+\alpha_2^2}Q_t(\vec
{0})\nonumber\\[-8pt]\\[-8pt]
&&\hspace*{39pt}
{}- 2C_t(\vec{\alpha})\e^{-(\alpha_1+\alpha_2)\sqrt{\lambda t}}
\e^{\alpha_1^2/2+\alpha_2^2/2} \bigr) =o(1).\nonumber
\end{eqnarray}
For these proofs, the explicit computations of mean and covariances\break
of~$z_t^{(1,2)}(\vec{a}_t)$ and $z_t^{(1,2)}=z_t^{ (1,2)}(1,1)$ are crucial.

To prove (\ref{aim-1}), we rewrite
%
%e3.90 ###
\begin{eqnarray}
\e^{-2\lambda t}M_t(\vec{\alpha}) &=&\e^{(a_1^s+a_2^s)\lambda t-2\lambda t-(\alpha_1+\alpha_2)\sqrt
{\lambda t}} \e^{-\lambda(a_1^s+a_2^s) t}\expec\bigl(z_t^{ (1,2)}(\vec{a})\bigr) \nonumber\\[-8pt]\\[-8pt]
&&{}-\e^{\alpha_1^2/2+\alpha_2^2/2}\e^{-2\lambda t}\expec\bigl(z_t^{(1,2)}(\vec{1})\bigr).\nonumber
\end{eqnarray}
Since $\vec{a}_t=(\e^{\alpha_1/\sqrt{s^2\lambda t}},\e^{\alpha
_2/\sqrt{s^2\lambda t}})
\to(1,1),$
%
%e3.91 ###
\begin{equation} \e^{-2\lambda t}\expec\bigl(z_t^{(1,2)}(\vec{1})\bigr)\to A,\qquad
\e^{-\lambda(a_1^s+a_2^s) t}\expec\bigl(z_t^{(1,2)}(\vec{a})\bigr) \to A,
\end{equation}
where $A$ is the limit in (\ref{asy-mean-two-vertex}) in
Lemma \ref{lemma:mean-two-vertex-charac}(a). By a second order
Taylor expansion,
%
%e3.92 ###
\begin{equation} \bigl(a_1^s \lambda t-\lambda t-\alpha_1\sqrt{\lambda
t}\bigr)= \lambda t\biggl(\e^{\alpha_1/\sqrt{\lambda t}}-1-\frac{\alpha
_1}{\sqrt{\lambda t}}\biggr) =\alpha_1^2/2+o(1).
\end{equation}
Together, these two asymptotics show that (\ref{aim-1}) holds.
The proof of (\ref{aim-2}) is identical, now using
Lemma \ref{lemma:mean-two-vertex-charac}(b) instead, and the fact that
the limit equals $A_{\operatorname{Cov}}(\vec{1}, \vec{1})$ for all contributions,
since $\vec{a}_t\rightarrow(1,1)$.
\end{pf}

\begin{pf*}{Completion of the proof of Theorem \ref{theo:two-vertex-charac-ratio}}
We can bound
%
%e3.93 ###
\begin{eqnarray}\\
 && \biggl|\frac{\sum_{i\in\SWG_t^{(1)}(x)} \sum_{j\in\SWG
_t^{ (2)}(y)}\chi_{ij}(t)}{\sum_{i\in\SWG_t^{(1)}} \sum_{j\in
\SWG_t^{(2)}}\chi_{ij}(t)} -\frac{\sum_{i\in\BP_t^{(1)}(x)} \sum
_{j\in\BP_t^{ (2)}(y)}\chi_{ij}(t)}{\sum_{i\in\BP_t^{(1)}} \sum
_{j\in\BP_t^{(2)}}\chi_{ij}(t)} \biggr|\nonumber\\
 &&\qquad  \leq2\frac{\e^{-2\lambda
t}\sum_{i\in\BP_t^{(1)}} \sum_{j\in\BP_t^{(2)}}\chi_{ij}(t)-\e
^{-2\lambda t}\sum_{i\in\SWG_t^{ (1)}} \sum_{j\in\SWG
_t^{(2)}}\chi_{ij}(t)}{\e^{-2\lambda t}\sum_{i\in\BP_t^{(1)}}
\sum_{j\in\BP_t^{(2)}}\chi_{ij}(t)}.\nonumber
\end{eqnarray}
The random variable in the denominator converges in probability
to\break
$\kappa W^{(1)} W^{(2)}>0$ by
Theorem \ref{theo:denom-as-convg}, so that it suffices to prove that
the numerator converges to 0 in
probability.

Denote
%
%e3.94 ###
\begin{equation} z_t^{n,{(1,2)}}=\sum_{i\in\SWG_t^{(1)}} \sum_{j\in
\SWG_t^{(2)}}\chi_{ij}(t).
\end{equation}
Then, similarly to Proposition \ref{prop:coupling-simul}(c), and
recalling that
$t_n = T_{12}\wedge t_n^*$ where $t_n^* = (2\lambda)^{-1}\log{n} +B$
for some $B>0$,
we obtain that $\sup_{t\leq t_n} (z_t^{n,(1,2)} - z_t^{{
(1,2)}})$ is tight.
From Theorem \ref{theo:wt-shwp} (whose proof has been completed since
it relies only on Theorem \ref{theo:denom-as-convg}, which was proved
in the previous section),
we know that the collision time $T_{12}$ is bounded by $t_n^*$ with
probability $1-o(1)$ as
$B\uparrow\infty$.
Therefore,
%
%e3.95 ###
\begin{eqnarray}
&& \e^{-2\lambda t}\sum_{i\in\BP_t^{(1)}} \sum_{j\in
\BP_t^{(2)}}\chi_{ij}(t)-\e^{-2\lambda t}\sum_{i\in\SWG_t^{
(1)}}\sum_{j\in\SWG_t^{(2)}}\chi_{ij}(t)\nonumber\\[-8pt]\\[-8pt]
&&\qquad  =\e^{-2\lambda
t}\bigl(z_t^{n,{(1,2)}} - z_t^{{(1,2)}}\bigr)\convp0.\nonumber
\end{eqnarray}
This completes the proof.
\end{pf*}
%t_n} (z_t^{n,(1,2)} - %z_t^{{(1,2)}})$?}

\begin{appendix}
\section*{Appendix: Auxiliary results}\label{sec-app-A}
\setcounter{Lemma}{0}
\setcounter{equation}{0}
In this section, we prove an auxiliary result on the asymptotics
of the measure $\tilde{\mu}_a(du)=\sum_{j=0}^{\infty} a^j\mu^{*j}(du)$.

\begin{Lemma}[(Asymptotics of density of $\tilde{\mu}_a$)]
\label{lem-dens-tildemu}
\begin{longlist}[(a)]
\item[(a)] Let
%
%e3.96 ###
\begin{equation} \e^{-\lambda a^s u}\sum_{j=1}^\infty a^j \mu
^{*j}(du) \equiv p_a(u) \,du.
\end{equation}
Then, for $u\in[0,1]$, there exists a constant $c$ such that
$p_a(u)\leq cu^{1/s-1}$, while, for $u\geq1$, $p_a(u)$ is bounded and
as $u\rightarrow\infty$,
%
%e3.97 ###
\begin{equation} \label{p(u)-conv} p_1(u)\rightarrow\lambda s.
\end{equation}
\item[(b)] The following scaling identity holds:
%
%e3.98 ###
\begin{equation} \label{pa-scaling} p_a(u)=a^{s}p_1(ua^s).
\end{equation}
In particular, when $a_u\rightarrow1$,
%
%e3.99 ###
\begin{equation} \label{pa-asy} p_{a_u}(u)\rightarrow\lambda s.
\end{equation}
\end{longlist}
\end{Lemma}

\begin{pf}
By (\ref{mu*j-form}), we have
%
%e3.100 ###
\begin{equation}
p_a(u)=\e^{-\lambda a^s u}\sum_{j=1}^{\infty} a^j\frac
{u^{j/s-1} \lambda^{j/s}}{\Gamma(j/s)}.
\end{equation}
This form immediately proves the identity in (\ref{pa-scaling}), and
therefore also
(\ref{pa-asy}) follows from (\ref{p(u)-conv}). Also, this form
immediately shows that $p_a(u)\leq cu^{1/s-1}$ for $u\in[0,1]$.
Thus, we are left to prove (\ref{p(u)-conv}).

By \cite{GraRyz65}, equation 8.327, we have that, as $z\rightarrow\infty$,
%
%e3.101 ###
\begin{equation} \label{Gamma-asy} z^{z-1/2}\e^{-z} \sqrt{2\pi}\leq
\Gamma(z)\leq z^{z-1/2}\e^{-z} \sqrt{2\pi} \biggl(1+\frac{1}{12 z}\biggr).
\end{equation}
Therefore,
%
%e3.102 ###
\begin{eqnarray}
p_1(u)&=&\bigl(1+o(1)\bigr)\e^{-\lambda u}\sum_{j=1}^{\infty}
\frac{1}{\sqrt{2\pi s/j}} u^{j/s-1} \lambda^{j/s} \e^{j/s}
(j/s)^{-j/s}\nonumber\\[-8pt]\\[-8pt]
&=&\bigl(1+o(1)\bigr)\lambda\e^{-\lambda u}(\lambda u)^{-1}\sum
_{j=1}^{\infty} \frac{\sqrt{j}}{\sqrt{2\pi s}} (\lambda u s \e
j^{-1})^{j/s},\nonumber
\end{eqnarray}
where the error term converges to 0 as $u\rightarrow\infty$.
Note that the right-hand side\ is a function of $\lambda u$, so that it suffices
to prove that
%
%e3.103 ###
\begin{equation} \label{aim-q} q(v)=\e^{-v}v^{-1} \sqrt{s} \sum
_{j=1}^{\infty} \frac{1}{\sqrt{2\pi j}} \e^{j/s \log{(v s \e
/j)}} \rightarrow s.
\end{equation}
For this, we note that $j\mapsto\e^{j/s \log{(v s \e/j)}}$ is
maximal when $j=sv$, where it takes the value $\e^{v}$. A second order
Taylor expansion shows that when $j-sv=x$, we have
%
%e3.104 ###
\begin{equation} \e^{j/s \log{(v s \e/j)}} =\e^{v} \e^{-x^2/(2s^2
v)}\bigl(1+o(1)\bigr).
\end{equation}
Performing the approximate Gaussian sum leads to the claim in (\ref{aim-q}).
\end{pf}
\end{appendix}

\section*{Acknowledgments}
Shankar Bhamidi would also like to thank the hospitality of \textsc{Eurandom}
where much of this work was done.
We would like to thank Jesse Goodman and Maren Eckhoff for many helpful
discussions.

% imsref loaded by dianan, 2011-02-15 09:53:35
%

\printaddresses

\end{document}